\documentclass[12pt]{amsart}
\usepackage{amscd,amsmath,amssymb,amsfonts}
\usepackage{hyperref}
\usepackage{color}

\def\timesover#1#2#3{\ \xymatrix@1@=0pt@M=0pt{ _{#1}&\times&_{#2} \\& ^{#3}&}\ }
\def\otimesover#1#2#3{\ \xymatrix@1@=0pt@M=0pt{ _{#1}&\otimes&_{#2} \\& ^{#3}&}\ }

\usepackage[all]{xy}
\theoremstyle{plain}
\newtheorem{thm}{Theorem}
\newtheorem{lem}[thm]{Lemma}
\newtheorem{cor}[thm]{Corollary}
\newtheorem{prop}[thm]{Proposition}

\theoremstyle{definition}
\newtheorem{defn}[thm]{Definition}

\newtheorem{rmk}[thm]{Remark}
\newtheorem{rmks}[thm]{Remarks}

\newtheorem{eig}[thm]{Properties}
\newtheorem{nota}[thm]{Notation}
\newtheorem{notas}[thm]{Notations}

\numberwithin{thm}{section}
\numberwithin{equation}{section}

\newcommand{\ml}[2]{\begin{multline}\label{#1}#2 \end{multline}}
\newcommand{\ga}[2]{\begin{gather}\label{#1}#2 \end{gather}}

\newcommand{\surj}{\twoheadrightarrow}
\newcommand{\inj}{\hookrightarrow}

\newcommand{\Hom}{{\rm Hom}}

\newcommand{\Spec}{{\rm Spec \,}}

\newcommand{\sC}{{\mathcal C}}
\newcommand{\sD}{{\mathcal D}}

\newcommand{\sO}{{\mathcal O}}

\newcommand{\sT}{{\mathcal T}}

\newcommand{\G}{{\mathbb G}}

\newcommand{\N}{{\mathbb N}}

\newcommand{\Q}{{\mathbb Q}}

\newcommand{\Rep}{{\rm Rep\hspace{0.1ex}}}
\newcommand{\id}{{\rm Id}}
\begin{document}

\title{The  fundamental groupoid scheme and applications}
\author{H\'el\`ene Esnault}
\address{
Universit\"at Duisburg-Essen, Mathematik, 45117 Essen, Germany}
\email{esnault@uni-due.de}
\author{Ph\`ung  H\^o Hai}
\address{
Universit\"at Duisburg-Essen, Mathematik, 45117 Essen, Germany
and Institute of Mathematics, Hanoi, Vietnam}
\email{hai.phung@uni-duisburg-essen.de}
\date{Nov. 15, 2006}
\thanks{Partially supported by  the DFG Leibniz Preis and   the DFG Heisenberg program}
\begin{abstract}
We define a linear structure on Grothendieck's arithmetic fundamental group $\pi_1(X, x)$  
of a scheme $X$ defined over a field $k$ of
characteristic 0. It allows us to link the existence of  sections of the Galois group ${\rm Gal}(\bar k/k)$ to $\pi_1(X, x)$  with the existence of a neutral fiber functor on the category which linearizes it. When applied to Grothendieck section conjecture, it allows us to find a $k$-structure on the universal covering, and $k$-rational pro-points at finite and infinite distance which lift given $k$-rational points. 
\end{abstract}
\maketitle
\begin{quote}

\end{quote}
\section{Introduction}
\noindent
For a connected scheme $X$ over a field $k$, Grothendieck defines in 
  \cite[Section~5]{Groth}   the arithmetic fundamental group as follows. He introduces the category ${\sf ECov}(X)$ of finite  \'etale coverings $\pi: Y\to X$, the Hom-Sets being $X$-morphisms. The choice of  a geometric point $x\in X$ 
defines a fiber functor $\pi \xrightarrow{\omega_x}  \pi^{-1}(x)$ with values in the category ${\sf FSets}$ of finite sets. He defines the arithmetic fundamental group with base point $x$  to be the automorphism group  $\pi_1(X,x)={\rm Aut}(\omega_x)$ of the fiber functor. It is an abstract group, endowed with the pro-finite topology stemming from its finite images in the permutation groups of $\pi^{-1}(x)$ as $\pi$ varies. The main theorem is the equivalence of categories ${\sf ECov}(X)\xrightarrow{\omega_x} \pi_1(X,x)\text{-\sf FSets}$ 
between the \'etale coverings and the finite sets acted on continuously by $\pi_1(X,x)$. 
This equivalence extends to pro-finite objects on both sides. When applied to the
set
$\pi_1(X, x)$, acted on by itself as a group via translations, it defines the universal pro-finite \'etale covering $\tilde{\pi}_x: \tilde{X}_x\to X$ based at $x$. It is a Galois covering of group $\pi_1(X, x)$ and  $\tilde{\pi}_x^{-1}(x)=\pi_1(X,x)$. Furthermore,
the embedding of the  sub-category ${\sf ECov}(k)$ of ${\sf ECov}(X)$,  consisting of the  \'etale coverings $X\times_k {\rm Spec}(L)\to X$ obtained from base change by a finite field extension $L\supset k$, with $L$ lying in the residue field $\kappa(x)$ of $x$, induces an augmentation map 
$\pi_1(X, x)\xrightarrow{\epsilon} {\rm Gal}(\bar{k}/k)$. Here $\bar{k}$ is the algebraic closure of $k$ in $\kappa(x)$. Then $\epsilon$
is surjective when $k$ is separably closed in $H^0(X, \sO_X)$. Grothendieck shows that the kernel of the augmentation is identified with $\pi_1(X\times_k \bar k, x)$, thus  one has an exact sequence
\ga{1.1}{1\to \pi_1(X\times_k \bar{k}, x)\to \pi_1(X, x)\xrightarrow{\epsilon}  {\rm Gal}(\bar{k}/k).}
Via the equivalence of categories, one has a factorization
\ga{1.2}{\xymatrix{ \tilde{X}_x \ar[r] \ar[dr]_{\tilde{\pi}_x} & X\times_k \bar{k} \ar[d]\\
& X
}
}
identifying the horizontal map with the pro-finite \'etale covering of $X\times_k \bar{k}$, based at $x$ viewed as a geometric point of $X\times_k \bar{k}$. In particular, $\tilde{X}_x$ is a $\bar{k}$-scheme. 
\\[.2cm]
The aim of this article is to linearize Grothendieck's construction in the case where $X$ is smooth and $k$ has characteristic 0. There is a standard way to do this using  local systems of $\Q$-vector spaces on $X_{{\rm \acute{e}t}}$, which is a $\Q$-linear abelian rigid tensor category (\cite[Section 10]{DeP}). 
Rather than doing this, we go to the $\sD$-module side of the Riemann-Hilbert correspondence, where we  have more flexibility for the choice of a fiber functor. 
Before describing more precisely our construction,  the theorems and some applications, let us first recall Nori's theory.\\[.2cm] 
For a proper, reduced scheme $X$ over a perfect field $k$, which is connected in the strong sense that $H^0(X, \sO_X)=k$, Nori defines in 
\cite[Chapter~II]{N2} the fundamental group scheme as follows. He introduces the category $\sC^N(X)$ of essentially finite bundles, 
which is the full sub-category of the coherent category on $X$, spanned by Weil-finite bundles $V$, that is for which there are
two polynomials $ f,g\in \N[T], f\neq g$, with the property that $f(V)$ is  isomorphic to $g(V)$. It is a
 $k$-linear  abelian rigid tensor category.
 The choice of a rational point $x\in X(k)$ defines a fiber functor $V \xrightarrow{\rho_x} V|_x$  with 
values in the category ${\sf Vec}_k$ of vector spaces over $k$, thus endows $\sC^N(X)$ with the structure of a neutral Tannaka category. Nori defines his fundamental group scheme with base point $x$ to be the  automorphism group 
 $\pi^N(X,x)={\rm Aut}^{\otimes}(\rho_x)$ of $\rho_x$.
 It is a $k$-affine group scheme, and its affine structure is pro-finite in the sense that its images in $GL(V|_x)$ are 0-dimensional (i.e. finite) group schemes over $k$ for all objects $V$ of $\sC^N(X)$. 
Tannaka duality 
(\cite[Theorem~2.11]{DeMil}) asserts the equivalence of categories $\sC^N(X)\xrightarrow{\rho_x} {\rm Rep} (\pi^N(X,x))$. The duality extends to the ind-categories on both sides. When applied to $k[\pi^N(X,x)]$ acted on by $\pi^N(X,x)$ via translations, it defines a principal bundle $\pi_{\rho_x}: X_{\rho_x}\to X$ under $\pi^N(X,x)$. One has $\pi_{\rho_x}^{-1}(x)=\pi^N(X,x)$. In particular,
the $k$-rational point $1 \in \pi^N(X,x)(k)$ is a lifting in $X_{\rho_x}$ of the $k$-rational point $x\in X(k)$ .
 On the other hand, one has the base change property for $\sC^N(X)$. Any $\bar{k}$-point $\bar{x}\to X$ lifting $x$
yields a base change isomorphism $\pi^N(X,x)\times_k \bar{k}\cong \pi^N(X\times_k \bar{k}, \bar{x})$, where  $k\subset \bar{k}$ is given by $\bar{x}\to x$. Denoting by  $\rho_{\bar{x}}$ the fiber functor
$V\xrightarrow{\rho_{\bar x}} V|_{\bar x}$ on $\sC^N(X\times_k \bar k)$, it implies an isomorphism 
$(X\times_k \bar{k})_{\rho_{\bar x}}\cong X_{\rho_x} \times_k \bar{k}$. 
Given those two parallel descriptions of $\tilde{\pi}_{x}: \tilde{X}_x\to X$ and  $\pi_{\rho_x}: X_{\rho_x}\to X$ in Grothendieck's and in Nori's theories, 
together with the base change property of $\sC^N(X)$,
one deduces $\pi_1(X\times_k \bar{k},x)\cong \pi^{\acute{e}t}(X\times_k \bar{k},x)(\bar{k})$, where $\pi^{\acute{e}t}(X,x)$ is the quotient of $\pi^N(X,x)$ obtained by considering the quotient pro-system of \'etale group schemes  (see \cite[Remarks~3.2, 2)]{EsHS}). So in particular, if $k$ has characteristic zero, one has  $\pi_1(X\times_k \bar{k},\bar x)\cong \pi^{N}(X\times_k \bar{k},\bar x)(\bar{k})$.
We conclude that $\pi_{\rho_{\bar x}}=\tilde{\pi}_{\bar x}$. 
So \eqref{1.2} has the factorization
\ga{1.3}{\xymatrix{ \ar[d] \tilde{X}_{\bar x} \ar[r] \ar[dr]_{\tilde{\pi}_{\bar x}} & X\times_k \bar{k} \ar[d]\\
 X_{\rho_x}\ar[r]_{\pi_{\rho_x}} & X
}
}
(See Theorem \ref{thm2.11} and Remarks \ref{rmk2.13}). \\[.2cm]
Summarizing, we see that Grothendieck's construction is very closed to a Tannaka construction, except that his category ${\sf ECov}(X)$ does not have a $k$-structure. On the other hand, the comparison of
$\pi_1(X,x)$ with $\pi^N(X,x)$, aside of the assumptions under which $\pi^N(X,x)$ is defined, that is $X$ proper, reduced and strongly connected, and $x\in X(k)$, requests the extra assumption that $k$   be algebraically closed. In particular, the sub-category
${\sf ECov}(k)$ of ${\sf ECov}(X)$ 
is not seen by Nori's construction, so the augmentation $\epsilon$ in \eqref{1.1} is not seen either, and there is no influence of the Galois group of $k$.  The construction only yields a $k$-structure on 
 $\pi_1(X\times_k \bar{k}, \bar x)$ under the assumption that $\bar x$ lifts a rational point. However, when this assumption is fulfilled, it yields a $k$-form of the universal pro-finite \'etale covering which carries a rational pro-point. 
\\[.2cm]
The purpose of this article is to reconcile the two viewpoints, using Deligne's more evolved Tannaka formalism as developed in \cite{DeGroth}. {F}rom now on, we will restrict ourselves to the characteristic 0 case. 
We  illustrate the first idea on the simplest possible example 
 $X=x={\rm Spec}(k)$.  Then certainly $\sC^N(X)$ is the trivial category, that is every object is isomorphic to a direct sum of the trivial object, thus $\pi^N(X,x)=\{1\}$. Let us fix $\bar{x}\to x$ corresponding to the choice of an algebraic closure $k\subset \bar{k}$. It defines  the fiber functor
$\sC^N(X)\xrightarrow{\rho_{\bar{x}}} {\sf Vec}_{\bar{k}}$ which assigns $V|_{\bar{x}}=V\otimes_k \bar{k} $ to $V$. It is a  non-neutral fiber functor. It defines a groupoid scheme ${\rm Aut}^{\otimes}(\rho_{\bar{x}})$  over $k$, acting transitively on $\bar{x}$ via $(t,s): 
{\rm Aut}^{\otimes}(\rho_{\bar{x}}) \to \bar{x}\times_k \bar{x}$. In fact 
${\rm Aut}^{\otimes}(\rho_{\bar{x}})=\bar{x}\times_k \bar{x}$ 
    and Deligne's Tannaka duality  theorem \cite[Th\'eor\`eme~1.12]{DeGroth} asserts that the trivial category $\sC^N(X)$ is equivalent to the representation category of the trivial $k$-groupoid scheme
$\bar{x}\times_k \bar{x}$ acting on $\bar{x}$.
We define $( \bar{x}\times_k \bar{x})_s$ to be  $\bar{x}\times_k \bar{x}$ viewed as a $\bar{k}$-scheme by means of the right projection $s$. 
Galois theory then implies that the $\bar{k}$-points of $( \bar{x}\times_k \bar{x})_s$ form a pro-finite group which is the Galois group of ${\rm Gal}(\bar{k}/k)$ where the embedding $k\subset \bar{k}$ corresponds to $\bar{x}\to x$. Thus the geometry here is trivial, but the arithmetic is saved via the use of a non-neutral fiber functor. This is the starting point of what we want to generalize. \\[.2cm]
In order to have a larger range of applicability to not necessarily proper schemes, we extend in section 2 Nori's definition to smooth schemes of finite type $X$ over $k$ which have the property that the field $k$, which we assume to be of characteristic 0, is exactly the field of constants of $X$, that is $k= H^0_{DR}(X)$. In particular, it forces the augmentation $\epsilon$ in \eqref{1.1} to be surjective.  We define the category ${\sf FConn}(X)$ to be the category of finite flat connections, that is of bundles $V$ with a flat connection $\nabla: V\to \Omega^1_X\otimes_{\sO_X} V$, such that $(V,\nabla)$ is a sub-quotient connection of a Weil-finite connection. The latter means 
that there are $f, g\in \N[T], f\neq g$ such  that $f((V,\nabla))$ is isomorphic to $g((V,\nabla))$.  The Hom-Sets are flat morphisms. ${\sf FConn}(X)$ is a $k$-linear rigid tensor category (see Definitions \ref{defn2.1} and \ref{defn2.3}). Theorem \ref{thm2.15} shows that 
if $X$ is smooth proper, with $k=H^0(X,\sO_X)$ of characteristic 0, then the forgetful functor ${\sf FConn}(X)\to \sC^N(X), \ (V,\nabla)\mapsto V$ is an equivalence of $k$-linear abelian rigid tensor categories. \\[.2cm]
A fiber functor $\rho: {\sf FConn}(X)\to {\sf QCoh}(S)$ in the quasi-coherent category  of a scheme $S$ over $k$ endows ${\sf FConn}(X)$ with a Tannaka structure. One defines $\Pi={\rm Aut}^{\otimes }(\rho)$ to be the  groupoid  scheme over $k$ acting transitively on $S$, with diagonal $S$-group scheme $\Pi^\Delta$. The fiber functor $\rho$ establishes an equivalence ${\sf FConn}(X)\xrightarrow{\cong \ \rho} {\rm Rep}(S: \Pi)$ (\cite[Th\'eor\`eme~1.12]{DeGroth}). Our first central theorem says that one can construct a universal covering in  this abstract setting, associated to $S$ and $\rho$,   as a direct generalization of the easily defined $\pi_{\rho_x}: X_{\rho_x}\to X$ 
sketched in \eqref{1.3}. 
\begin{thm}[See precise statement in Theorem \ref{thm2.7}] \label{thm1.1}
 There is a diagram of $k$-schemes
\ga{1.4}{\xymatrix{& \ar[ddl]_{p_\rho} X_{\rho} \ar[d]^{s_\rho} \ar[ddr]^{\pi_\rho}\\
& S\times_k X \ar[dl]^{p_1}  \ar[dr]_{p_2}\\
S &  & X}}
with the following properties.
\begin{itemize}
\item[1)]$s_\rho$ is a 
$\Pi^\Delta$-principal bundle, that is
$$X_\rho \times_{S\times_k X}X_\rho\cong \Pi^\Delta \times_S X_\rho.$$
\item[2)] $R^0 (p_\rho)_{* DR}(X_\rho/S):= R^0(p_\rho)_*(\Omega^\bullet_{X_\rho/S})=\sO_S$. 
\item[3)] For all objects $N=(W,\nabla)$ in ${\sf FConn}(X)$, the connection
$$\pi_\rho^*N: \pi_\rho^*W\to 
 \pi_\rho^*W\otimes_{\sO_{X_\rho}} \Omega^1_{X_\rho/S} $$
relative to $S$  is endowed with a functorial isomorphism with the relative connection $$p_\rho^*\rho(N)=
(p_\rho^{-1}\rho(N)\otimes_{p_\rho^{-1}\sO_S} \sO_{X_\rho}, 1\otimes d)$$ 
which is generated by the relative flat sections $p_\rho^{-1}\rho(N)$.
\item[4)] One recovers the fiber functor $\rho$ via $X_\rho$ by an isomorphism
$\rho(N)\cong R^0(p_\rho)_{*DR}(X_\rho/S, \pi_\rho^*N)$ which  is compatible with all morphisms in ${\sf FConn}(X)$. In particular, the data in \eqref{2.17} are equivalent to the datum $\rho$ (which defines $\Pi$). 

\end{itemize}
\end{thm}
\noindent
Furthermore this construction is functorial in $S$. While applied to $S=X, \ \rho((V,\nabla))=V$, $X_\rho$ is nothing but $\Pi$ with $(p_\rho, \pi_\rho)=(t,s)$ (Definition \ref{defn2.8}, 2)).  While applied to $S={\rm Spec}(\bar k), \ \rho((V,\nabla)):=\rho_{\bar x}((V,\nabla))=V|_{\bar{x}}$ for a geometric point $\bar{x}\to X$ with residue field $\bar{k}$, then $\pi_{\rho_{\bar x}}$ is an \'etale pro-finite covering, which will turn out to be Grothendieck's universal pro-finite \'etale covering. Let us be more precise. \\[.2cm]
Theorem \ref{thm1.1} is proven by showing  that $s_*\sO_\Pi$ is a representation of the groupoid scheme $\Pi$, by analyzing some of its properties, and by  translating them  via Tannaka duality. It invites one to consider the scheme $\Pi$ considered as a $S$-scheme via $s$, which we denote by $\Pi_s$. In case $S={\rm Spec}(\bar{k}), \ \rho=\rho_{\bar x}$, we denote 
${\rm Aut}^{\otimes}(\rho_{\bar x})$ by $\Pi(X,\bar x)$. We show 
\begin{thm}[See Theorem  \ref{thm4.4}] \label{thm1.2}
The rational points $\Pi(X,{\bar x})_s(\bar k)$ carry a structure of a pro-finite group. One has 
an exact sequence of pro-finite groups
\ga{1.5}{1\to \Pi(\bar X,{\bar x})(\bar{k})\to \Pi(X,{\bar x})_s(\bar k) 
\to ({\rm Spec}(\bar k)\times_k {\rm Spec}(\bar k))_s(\bar{k})
\to 1 }
together with an identification of exact sequences of pro-finite groups 
\ga{1.6}{\xymatrix{ 1\ar[r] & \Pi(\bar X,\bar x)(\bar k)\ar[r] \ar[d]_{ =} & 
 \Pi(X,\bar x)_s(\bar k) \ar[d]_{ = }\ar[r] &
({\rm Spec}(\bar k)\times_k {\rm Spec}(\bar k))_s(\bar{k}) 
\ar[r] \ar[d]_{= }& 1\\
1\ar[r] & \pi_1(\bar X,\bar x) \ar[r]
& \pi_1(X,\bar x)\ar[r]^{\epsilon} & {\rm Gal}(\bar{k}/k) \ar[r] &1
}}\end{thm}
\noindent
A corollary is that $\pi_{\rho_{\bar x}}: X_{\rho_{\bar x}}\to X$ is precisely the universal pro-finite \'etale covering of $X$ based at $\bar x$ (see Corollary \ref{cor4.5}). 
One interesting consequence of Theorem \ref{thm1.2} is that it yields a Tannaka interpretation of the existence of a splitting of the augmentation $\epsilon$. Indeed, \eqref{1.5} does not have to do with the choice of the category ${\sf FConn}(X)$. More generally, if $(t,s): \Pi \to {\rm Spec}(\bar k)\times_k {\rm Spec}(\bar k) $ is a  $k$-groupoid 
scheme acting transitively upon ${\rm Spec}(\bar k)$, then we show 
\begin{thm}[See Theorem \ref{thm3.2}] \label{thm1.3}
\begin{itemize}
\item[1)] There exists a group structure on $\Pi_s(\bar k)$
such that the map 
$$(t,s)|_{\Pi_s}: \Pi_s(\bar k)\to ({\rm Spec}(\bar{k})\times_k {\rm Spec}(\bar{k}))_s(\bar k)\cong 
{\rm Gal}(\bar{k}/k)$$
 is a group homomorphism. 
\item[2)]Splittings of $(t,s)|_{\Pi_s}: \Pi_s(\bar k)\to ({\rm Spec}(\bar{k})\times_k 
{\rm Spec}(\bar{k}))_s(\bar k)\cong {\rm Gal}(\bar{k}/k)$ 
as group homomorphisms are in one to one correspondence  with splittings of 
$(t,s): \Pi\to {\rm Spec}(\bar{k})\times_k {\rm Spec}(\bar{k})$ as $k$-affine groupoid 
scheme homomorphisms. 
\item[3)]  The projection $(t,s): \Pi\to {\rm Spec}(\bar{k})\times_k {\rm Spec}(\bar{k})$ has a section of  groupoid schemes over 
$k$ if and only if $\sC$ has a neutral fiber functor.
More precisely,
there is a 1-1 correspondence between neutral fiber functors of $\sC$ and splittings of $(t,s)$
up to an inner conjugation of $\Pi$ given by an element of $\Pi^\Delta(\bar k)$. 
\end{itemize}
\end{thm}
\noindent
In particular we show that splititngs of $\epsilon$ in \eqref{1.6} up to conjugation with 
$\pi_1(\bar X, \bar x)$ are in one to one correspondence with neutral fiber functors of ${\sf FConn}(X)\to {\sf Vec}_k$ (see Corollary \ref{cor4.6}). \\[.2cm]
We now discuss applications. Grothendieck, in his letter to G. Faltings  dated June 27-th, 1983 \cite{GroFa}, initiated a program to recognize hyperbolic algebraic curves defined over function fields $k$ over $\Q$ via the exact sequence \eqref{1.1}. His  anabelian conjectures have essentially been proven, with the notable exception of the so-called section conjecture. It predicts that if $X$ is a smooth projective curve of genus $\ge 2$ defined over $k$ of finite type over $\Q$, then conjugacy classes of sections of $\epsilon$ in \eqref{1.1} are in one to one correspondence with 
rational points of $X$. Our theorems \ref{1.2} and \ref{1.3} have then the following consequence.
\begin{thm}[See Theorem \ref{thm5.1}] \label{thm1.4}
If Grothendieck's section conjecture is true, then the  set of $k$-rational points of a  smooth projective curve $X$ 
of genus $g\ge 2$ over $k$ of finite type over $\Q$  
is in bijection with the set of neutral fiber functors ${\sf FConn}(X)\to {\sf Vec}_k$. 
\end{thm}
\noindent
Grothendieck's prediction is more precise. It says in particular that not only a section of $\epsilon$ should come from a  rational point, but also that the universal pro-finite \'etale covering of $\tilde{X}_{\bar{x}}$ based at $\bar{x}\in X(\bar{k})$, which is a $\bar{k}$-scheme, has a $k$-form which carries a rational point. The $k$-linear   structure on Grothendieck's arithmetic fundamental group allows us to 
answer positively  the easier part of this question.
\begin{prop}[See precise statement in  Proposition \ref{prop5.3}] \label{prop1.5} 
If $\bar{x}$ lies above $x\in X(k)$, then there is a $k$-scheme $X_x$ with $X_x\times_k \bar{k}=\tilde{X}_{\bar{x}}$ and with a $k$-rational pro-point above $x$. For any other geometric point $\bar{y}\to X$ with residue field $\bar{k}$, $X_x$ yields a $k$-form of $\tilde{X}_{\bar{y}}$, and thus a $k$-rational porpoint in it. 

\end{prop}\noindent
By Faltings' theorem, if $X'$ is a smooth projective curve of genus $\ge 2$ over a number field $k$,  which is large enough such that $X'(k)\neq \emptyset$, then $X'(k)$ is finite. So setting $X=X'\setminus X(k)$, one has of course $X(k)=\emptyset$, and there are tangential base points at $\infty=X'\setminus X$. They yield neutral fiber functors on ${\sf FConn}(X)$ as defined by Deligne \cite[Section 15]{DeP}, and Katz
\cite[Section~2.4]{KaGal}, thus sections of $\epsilon$. This was of course understood without the $k$-linear structure we defined on $\pi_1(X,x)$. But in light of the abstract construction of the universal  covering 
in Theorem \ref{thm1.1}, which, when 
 the   fiber functor is neutral, yields a cartesian square
\ga{1.7}{
\xymatrix{ X_{\rho\times_k \bar{k}}\cong \bar{k}\times_k X_\rho \ar[d] \ar[r] 
\ar@{}[rd]|\Box
& \bar{k}\times_k X \ar[d] \\
X_\rho  \ar[r] &  X}
}
(see Theorem \ref{thm2.11}), we discuss 
Grothendieck's section conjecture in the non-proper case 
for sections arising this way. 
\begin{thm}[See precise statement in Theorem \ref{thm5.7}] \label{thm1.6}
Let $\bar{y}\to X$ be a geometric point with residue field $\bar{k}$. Then 
$\tilde{X}_{\bar{y}}$ has a $k$-structure $X_\rho$ which has the property that the normalization of $X'$ in the function field of $X_\rho$ has a $k$-rational pro-point lifting $x$. 
\end{thm}
\noindent
This yields  a positive answer to the conjecture formulated in \cite{GroFa}, page 8,  by
Grothendieck. So we can say that Tannaka methods, which are purely algebraic, and  involve neither the geometry of $X$ nor the arithmetic of $k$, can't detect rational points on $X$. But they allow one, once one has a rational point on $X$, to follow it in the prosystem which underlies the Tannaka structure, in particular on Grothendieck's universal covering. 
\\[.2cm]
They also allow to generalize in a natural way Grothendieck's conjecture 
to general fiber functors $\rho: {\sf FConn}(X)\to {\sf QCoh}(S)$. (See Remark \ref{rmk5.8}, 3) and Remark \ref{rmk2.13} 2).). Theorem \ref{thm1.1} 2) says that $\rho$ is always cohomological in the sense that 
there is an isomorphic fiber functor which is cohomologically defined. 
  Theorem \ref{thm1.3} suggests to ask under which conditions on the geometry of $X$ and the arithmetic of $k$ 
one may say that $X$ has a $S$-point if and only if an $S$-valued $\rho$ exists. Further Theorem \ref{thm1.1} allows one to ask for a generalization to $S$-pro-points
of Theorem \ref{thm1.6}.  
\\[.2cm]
{\it Acknowledgements:} We thank Alexander Beilinson and Gerd Faltings 
for clarifying discussions,  for their interest and their encouragements. 

\section{Finite connections}
In this section, we introduce the category of finite connections in the spirit of Weil  and 
follow essentially verbatim Nori's developments in \cite[I, Section~2.3]{N2} for the main properties.  

\begin{defn} \label{defn2.1}
Let $X$ be a smooth scheme of finite type over a field $k$ of characteristic 0, with the property that $k=H^0_{DR}(X):=H^0(X, \sO_X)^{d=0}$. 
The category ${\sf Conn}(X)$ of flat connections has objects $M=(V, \nabla)$ where $V$ is a vector bundle (i.e. a locally free coherent sheaf) and $\nabla: V\to \Omega^1_X\otimes_{\sO_X}  V$ is a flat connection. The Hom-Sets are flat morphisms $f: V\to V'$, i.e. morphisms of coherent sheaves which commute with the connection. The rank of $M$ is the rank of the underlying vector bundle $V$. \end{defn}
Throughout this section,  we  fix a scheme $X$ as in Definition
\ref{defn2.1},
except in  Proposition \ref{prop2.14} and in Theorem \ref{thm2.15}, 
where we  relax the smoothness condition, but request the scheme to be proper. 
\begin{eig} \label{eig2.2} \begin{itemize}
\item[1)]Standard addition, tensor product and exact sequences of  connections endow
${\sf Conn}(X)$ with the structure of an abelian  $k$-linear rigid tensor category. 
\item[2)] ${\sf Conn}(X)$ is locally finite, that is the Hom-Set of any two objects is a finite dimensional vector space
over $k$ and each object has a decomposition series of finite length.
\item[3)] Every object $M$ in ${\sf Conn}(X)$ has only 
finitely many non-isomorphic sub-objects and quotient objects. \end{itemize}
\end{eig}
\begin{proof}
1) and 2) are classical (see e.g. \cite{KaGal}). \\[.2cm]
As for 3) we argue by induction on the rank. The only non-zero sub-connection of a rank 1 connection $M$ is $M$ itself. Assume rank $M>1$. Assume there is a sub-object $N\subset M$, thus rank $N$ and rank $M/N$ are $<$ rank $M$. Thus any sub-object $P\subset M$ is an extension of ${\rm Im} P\subset M/N$ by $P\cap N$, and we are reduce to see that for any objects $M_i$, ${\rm Ext}^1_{{\sf Conn}(X)} (M_1, M_2)$ is a finite dimensional $k$-vector space. But ${\rm Ext}^1$ is computed by de Rham cohomology ${\rm Ext}^1_{{\sf Conn}(X)} (M_1, M_2)=H^1_{DR}(X, M_1^\vee \otimes M_2)$ (see e.g \cite[ Proposition~2.2]{EP}) and de Rham cohomology is finite dimensional over $k$. 
\end{proof}
\begin{defn} \label{defn2.3}
${\sf FConn}(X)$ is defined to be the full sub-category of ${\sf Conn}(X)$ 
spanned by Weil-finite objects, where $M$ is Weil-finite if there are polynomials 
$f, g\in \N[T], f\neq g$, such that $f(M)$ is isomorphic to $g(M)$, where 
$T^n(M):=M\otimes \ldots \otimes M$ $n$-times and $nT(M):=M\oplus \ldots \oplus M$ 
$n$-times. Thus an object in ${\sf FConn}(X)$ is a sub-quotient in 
${\sf Conn}(X)$ of a Weil-finite object. 
\end{defn}
\begin{eig} \label{eig2.4}
\begin{itemize}
\item[1)] An object $M$ of ${\sf Conn}(X)$ is Weil-finite if and only if $M^\vee$ is Weil-finite. 
\item[2)] An object $M$ in ${\sf Conn}(X)$  is Weil-finite  if and only 
the collection of all isomorphism classes of indecomposable objects in all the tensor powers $M^{\otimes n}, n\in \N\setminus \{0\}$, is finite.
\item[3)] If an object $M$ in ${\sf Conn}(X)$ lies in ${\sf FConn}(X)$, then 
 there are finitely many isomorphism classes of indecomposable,
  in particular of simple, objects  in the tensor full sub-category 
  $\langle M\rangle$ of $  {\sf FConn}(X)$ 
 generated by $M$.
\item[4)] Let $L\supset k$ be a finite field extension, we shall denote
$X\times_k\Spec(L)$ by $X\times_kL$ for short. Denote by $\alpha: X\times_kL\to X$
the base change morphism.  Then $\alpha^* \alpha_* M\surj M$, $M\subset \alpha_*\alpha^*M$ for every object in ${\sf Conn}(X\times_k L), {\sf FConn}(X\times_k L)$ resp. ${\sf Conn}(X), {\sf FConn}(X)$.
\end{itemize}
\end{eig}
\begin{defn} \label{defn2.5} A locally finite abelian $k$-linear rigid tensor category which has finitely many isomorphism classes of indecomposable objects is called finite. \end{defn}

\begin{proof}[Proof of Properties \ref{eig2.4}]
For 1), $f(M)\cong g(M)$ if and only if $f(M^\vee)\cong g(M^\vee)$. 
For 2), one argues as in 
 \cite[I, Lemma~3.1]{N2}. One introduces the naive $K$ ring $K^{{\rm Conn}}(X)$ which is  spanned by 
 isomorphism classes of objects $M$ of ${\sf Conn}(X)$ modulo the relation $[M]\cdot [M']=[M\otimes M']$ and $[M]+[M']=[M\oplus M']$.
  Since in ${\sf Conn}(X)$ every object has a decomposition in a direct sum of finitely many indecomposable objects, 
  which is unique up to isomorphism, $K^{\sf Conn}(X)$ is freely spanned by classes $[M]$ with $M$ indecomposable. Then the proof of 2) goes word by word as in loc. cit. 
 Properties \ref{eig2.2} 2).
 Now  2) implies 3) and  4) is obvious.
\end{proof}

\noindent We will show in the sequel the converse of 3) using Tannaka duality. 
\\[.2cm]
We first recall some  basic concepts of groupoid schemes. Our ongoing reference
is \cite{DeGroth}.\\[.2cm]
An affine groupoid scheme $\Pi$ over $k$ acting transitively over $S$
is a $k$-scheme equipped with a faithfully flat morphism $(t,s):\Pi\longrightarrow
S\times_kS$ and 
\begin{itemize}\item [-] the multiplication $m:\Pi\timesover stS \Pi\longrightarrow \Pi$,
which is a morphism of $S\times_kS$-schemes,
\item [-] the unit element morphism $e:S\longrightarrow \Pi$,
which is a morphism of $S\times_kS$-schemes,
where $S$ is considered as an $S\times_kS$-scheme by means of the diagonal morphism,
\item[-] the inverse element morphism $\iota:\Pi\longrightarrow \Pi$ of $S\times_kS$-schemes,
which interchanges the morphisms $t,s$: $t\circ \iota=s$, $s\circ\iota=t$,\end{itemize}
satisfying the following conditions:
\ga{2.1}{m(m\times\id)=m(\id\times m),\\ \notag 
m(e\times\id)=m(\id\times e)=\id,\\
\notag m(\id\times \iota)=e\circ t, \quad m(\iota\times \id)=e\circ s.
}
One defines the diagonal group scheme $\Pi^\Delta$ over $S$ by the cartesian diagram
\ga{2.2}{\xymatrix{ \Pi^\Delta \ar[r] \ar[d] \ar@{}[rd]|\Box& \Pi \ar[d]^{(t,s)} \\
S \ar[r]_{{\Delta}\quad}&  S\times_k S}} 
where $\Delta: S\to S\times_k S$ is the diagonal embedding, and  the $S$-schemes $\Pi_s, \Pi_t$ by
\ga{2.3}{\Pi_s=S-{\rm scheme} \ \Pi\xrightarrow{s} S, \
\Pi_t=S-{\rm scheme} \ \Pi\xrightarrow{t} S .
}
For $(b,a): T\to S\times_k S$, one defines the $T$-scheme   $\Pi_{b,a}$
by the cartesian diagram
\ga{2.4}{\xymatrix{\Pi_{b,a}\ar[d] \ar[r] \ar@{}[rd]|\Box
 & \Pi \ar[d]^{(t,s)} \\
 T\ar[r]_{(b,a)\quad}&  S\times_k S }}
For three $T$ points $a,b,c: T\to S$ of $S$,  the  multiplication morphism $m$ yields the map
\ga{2.5}{\Pi_{c,b}\times_T\Pi_{b,a} 
\to \Pi_{c,a}
.}
By definition, the (abstract) groupoid $\Pi(T)$ has for objects morphisms $a:T\to S$ and
 for morphisms between $a$ and $b$ the set $\Pi_{b,a}(T)$ of $T$-points
 of $\Pi_{b,a}$. The composition law is given by \eqref{2.5}.\\[.2cm]
Moreover, the  multiplication $m$ induces an $S\times_k S$-morphism
\ga{}{ \mu: \Pi^\Delta \timesover{}{}S \Pi_t =\Pi^\Delta \timesover { }tS
\Pi
  \to \Pi\timesover{}{}{S\times_k S} \Pi
 \notag} 
by the following rule: for 
\ga{}{\xymatrix{\ar[dr]_b T\ar[r]^f & \Pi \ar[d]^t\\
&  S
} \qquad\qquad
\xymatrix{\ar[dr]_b T\ar[r]^g & \Pi^\Delta \ar[d]\\
&  S
}\notag
}
with $s\circ f=a : T\to S$, one has $f\in \Pi_{ba}(T),\  g\in \Pi_{bb}(T)$, thus one can compose 
$gf\in \Pi_{ba}(T)$ (as morphisms in the groupoid $\Pi(T)$) and  define
\ga{}{\mu(g,f)=(gf, f).\notag
}
One has 
\ga{2.6}{\mu \ {\rm defines \ a \ principal \ bundle \ structure} \\
{\rm i.e.} \ \Pi^\Delta \times_S \Pi_t
\cong \Pi\times_{S\times_k S} \Pi\
.\notag}
Indeed one has to check that $\mu$ is an isomorphism which is a local property on $S$.
In  \cite[Lemma~6.5]{EP} this fact was checked for the case $S=\Spec(K),$ $K\supset k$
a field, the proof can be easily extended for any $k$-algebras.
\\[.2cm]
We notice that those notations apply in particular to $\Pi=S\times_k S$ and one has
\ga{2.7}{{\sf Vec}_{k}={\rm Rep}(S: S\times_k S).}
In this case, $(S\times_kS)^\Delta=S$ via the diagonal embedding,
 and $(S\times_k S)_s$ is the scheme $S\times_kS$ viewed as a $S$-scheme via the second projection. 
\\[.2cm]
Let $V$ be a quasi-coherent sheaf on $S$. A representation $\chi$ of $\Pi$ on $V$ is a family of maps
\ga{2.8}{\chi_{b,a}:\Pi_{b,a}(T)\longrightarrow {\rm Iso}_T( a^*V\to b^*V)
 }
 for each $(b,a):T\to S\times_kS$,
which is compatible with the multiplication \eqref{2.7} and the base change $T'\to T$. 
 By definition, $\Pi_{t,s}$ is the fiber product 
$\Pi\times_{S\times_k S} \Pi$, seen as a $\Pi$-scheme via the first projection. Thus 
the identity morphism $\id_\Pi:\Pi\to \Pi$ can be seen as the morphism
between the two objects $t,s:\Pi\to S$ in the groupoid $\Pi_{t,s}(\Pi).$ This is the universal morphism.
 Then  $\chi$ is determined by the image of $\id_\Pi$ under $\chi_{t,s}$
  $$\chi_{t,s}(\id_\Pi): s^*V\stackrel\cong\longrightarrow t^*V.$$ 
The category of quasi-coherent representations of $\Pi$ on $S$  is denoted by $\Rep(S:\Pi)$.
Each representation of $\Pi$ is the union of its finite rank representations. That is
$\Rep(S:\Pi)$ is the ind-category of the category $\Rep_f(S:\Pi)$ of representations
on coherent sheaves on $S$. The category $\Rep_f(S:\Pi)$ is an abelian rigid tensor category. Its unit
object is the trivial representation of $\Pi$ on $\sO_S$.
On  the other hand, since $\Pi$ is faithfully flat over $S\times_k S$, one has 
\ga{2.9}{{\rm End}_{\Rep(S:\Pi)}(S)=k.}
Thus $\Rep_f(S:\Pi)$ is a $k$-linear abelian rigid tensor category.
\\[.2cm]
Thinking of $S\times_k S$ as a $k$-groupoid scheme acting on $S$, one has the projections $t$
 and $s$. In order not to confuse notations, we  denote the left projection
 $S\times_k S\to S$ by $p_1$ and by $p_2$ the right one. Thus $p_1=t, \ p_2=s$.  
\begin{thm}\label{thm2.6} Let $\Pi$ be an affine groupoid scheme over $k$ acting 
transitively on a $k$-scheme $S$. 
\begin{itemize}
\item[1)]
The formula \eqref{2.7} defines ${p_2}_*\sO_{S\times_k S}$ and 
$s_*\sO_{\Pi}$ as  objects in ${\rm Rep}(S: \Pi)$.
\item[2)]  ${p_2}_*\sO_{S\times_k S}$ and $s_*\sO_{\Pi}$ are algebra objects 
in ${\rm Rep}(S: \Pi)$.
\item[3)] The inclusion ${p_2}_*\sO_{S\times_k S} \to s_*\sO_{\Pi}$
of quasi-coherent sheaves on $S$ yields a ${p_2}_*\sO_{S\times_k S}$-algebra structure on $s_*\sO_{\Pi}$.
\item[4)] The maximal trivial sub-object of $s_*\sO_{\Pi}$ in ${\rm Rep}(S: \Pi)$ is the sub-object ${p_2}_*\sO_{S\times_k S}$.
\end{itemize}
\end{thm}
\begin{proof}
We first prove 1). Fixing $f\in \Pi_{b,a}(T)$, \eqref{2.7} yields a morphism $\Pi_{bc}\to \Pi_{ac}$ and consequently 
an $\sO_T$-algebra homomorphism $a^*s_*\sO_\Pi\to b^*s_*\sO_\Pi$. It is easily checked 
that this yields a representation of $\Pi$ on $s_*\sO_{\Pi}$. 
In particular,  the groupoid $S\times_kS$ has a
representation in ${p_2}_*\sO_{S\times_kS}$. The homomorphism of groupoid
schemes $(t,s):\Pi\longrightarrow S\times_kS$ yields a representation of
$\Pi$ in ${p_2}_*\sO_{S\times_kS}$, which is in fact trivial.
For 2), the algebra structure stemming from the structure 
sheaves is compatible with the $\Pi$-action. Then 3) is trivial. \\[.2cm]
As for 4), one 
notices that this property is local with respect to $S$, thus we can assume that
$S=\Spec R$, where $R$ is a $k$-algebra. Then $\sO_{\Pi} $ is
a $k$-Hopf algebroid acting on $R$. One has an isomorphism of categories
\ga{2.10}{\Rep(S:\Pi)\cong{\sf Comod}(R, \sO_{\Pi})}
where ${\sf Comod}(R,\sO_{\Pi})$ denotes the category of right $\sO_\Pi$-comodules in ${\sf Mod}_R$.
Under this isomorphism, the representation of $\Pi$ in $s_*\sO_{\Pi}$ corresponds
to the coaction of $\sO_{\Pi}$ on itself by means of the coproduct $\Delta:\sO_{\Pi}
\longrightarrow \sO_{\Pi} \otimesover stR \sO_{\Pi}$.
For $M\in \Rep(S:\Pi)$, one has the  natural isomorphism
\ga{2.11}{\Hom_{\Rep(S:\Pi)}(M,s_*\sO_{\Pi})\cong \Hom_{R}(M,R), \ f\mapsto e\circ f}
where $e:\sO_\Pi \to R $ is the unit element morphism. Indeed, the inverse to this map is 
given by the coaction of 
$\delta:M\to M\otimesover{}tR \sO_\Pi$ as follows
\ga{2.12}{\big(\varphi:M\to R\big) \mapsto \big((\varphi\otimes\id)\delta: 
M\to M\otimes_t\sO_\Pi\to \sO_\Pi\big).}
In particular, if we take $M=R$, the trivial representation of $\Pi$, then
we see that
\ga{2.13}{
 \Hom_{{\sf Comod}(R: \sO_\Pi)}(R, \sO_\Pi)\cong R. }
Notice that in ${\sf Comod}(R,\sO_\Pi)$, one has
$\Hom_{{\sf Comod}(R,\sO_\Pi)}(R,R)\cong k.$
We conclude that the maximal trivial sub-comodule of $\sO_\Pi$ is $R\otimes_kR$.
This finishes the proof. \end{proof}

\noindent In order to apply the theory of Tannaka duality as  developed by 
Deligne in \cite{DeGroth}, one needs a fiber functor. 
Let $S$ be a scheme defined over $k$ and let
\ga{2.14}{\rho: {\sf FConn}(X)\to {\sf QCoh}(S)}
be a fiber functor with values in the category of quasi-coherent sheaves over $S$, 
as in \cite[Section~1]{DeGroth}. Such a $\rho$ exists. For example, a very 
tautological one is provided by
$S=X,$ $\rho=\tau$ (for tautological): 
$$ {\sf FConn}(X)\to {\sf QCoh}(X), (V,\nabla)\mapsto V.$$
Then the functor from $S\times_k S$-schemes to ${\sf Sets}$, 
which assigns to any $(b,a): T\to S\times_k S$ the set of natural isomorphisms 
${\sf Iso}_T^{\otimes}(a^*\rho, b^*\rho)$ of fiber functors to
${\sf QCoh}(T)$, is representable by the affine groupoid scheme 
$$(t,s): \Pi:={\rm Aut}^{\otimes }(\rho) \to S\times_k S$$ defined over 
$k$ and  acting transitively over $S$ (\cite[Th\'eor\`eme~1.12]{DeGroth}):
\ga{2.15}{{\sf Iso}_T^{\otimes}(a^*\rho,b^*\rho)\cong {\sf Mor}_{S\times_kS}(T,\Pi).}
Tannaka duality asserts in particular that $\rho$ induces an equivalence between
${\sf FConn}(X)$ and the category $\Rep_f(S:\Pi)$ of representations of $\Pi$
in coherent sheaves over $S$.
This equivalence extends to an equivalence between the ind-category
$\text{\sf Ind-FConn}(X)$ of connections which are union of finite sub-connections
and the category $\Rep(S:\Pi)$ of representations of $\Pi$ in quasi-coherent
sheaves on $S$.\\[.2cm]
It follows from \eqref{2.15} for $T=S$, $(b,a)=(t,s)$ that 
there is an isomorphism 
\ga{2.16}{s^*\rho\cong t^*\rho.}
We now translate via Tannaka duality the assertions of Theorem \ref{thm2.6}.
\begin{thm} \label{thm2.7}
Let $X$ be a smooth scheme of finite type defined over a field $k$ of characteristic 0.
 Let $S$ be a $k$-scheme and $\rho: {\sf FConn}(X)\to {\sf QCoh}(S)$ be a fiber functor.
Let $\Pi =\rm{Aut}^\otimes(\rho)$  be the corresponding Tannaka $k$-groupoid scheme acting on $S$. 
 Then there is a diagram of $k$-schemes
\ga{2.17}{\xymatrix{& \ar[ddl]_{p_\rho} X_{\rho} \ar[d]^{s_\rho} \ar[ddr]^{\pi_\rho}\\
& S\times_k X \ar[dl]^{p_1}  \ar[dr]_{p_2}\\
S &  & X}}
with the following properties.
\begin{itemize}
\item[1)]$s_\rho$ is a 
$\Pi^\Delta$-principal bundle, that is
$$X_\rho \times_{S\times_k X}X_\rho\cong \Pi^\Delta \times_S X_\rho .$$
\item[2)] $R^0 (p_\rho)_{* DR}(X_\rho/S):= R^0(p_\rho)_*(\Omega^\bullet_{X_\rho/S})=\sO_S$. 
\item[3)] For all objects $N=(W,\nabla)$ in ${\sf FConn}(X)$, the connection
$$\pi_\rho^*N: \pi_\rho^*W\to 
 \pi_\rho^*W\otimes_{\sO_{X_\rho}} \Omega^1_{X_\rho/S} $$
relative to $S$  is endowed with an  isomorphism with the relative connection  $$p_\rho^*\rho(N)=
(p_\rho^{-1}\rho(N)\otimes_{p_\rho^{-1}\sO_S} \sO_{X_\rho}, \id\otimes d_{X_\rho/S})$$ 
which is generated by the relative flat sections $p_\rho^{-1}\rho(N)$. 
This isomorphism is compatible with all morphisms in ${\sf FConn}(X)$.
\item[4)] One recovers the fiber functor $\rho$ via $X_\rho$ by an isomorphism
$\rho(N)\cong R^0(p_\rho)_{*DR}(X_\rho/S, \pi_\rho^*N)$ which  is compatible with all morphisms in ${\sf FConn}(X)$. In particular, the data in \eqref{2.17} are equivalent to the datum $\rho$ (which defines $\Pi$). 
\item[5)] This construction is compatible with base change: if $u: T\to 
S$ is a morphism of $k$-schemes, and $u^*\rho: {\sf FConn}(X)\to {\sf QCoh}(T)$ is the composite fiber functor, then one has a cartesian diagram
\ga{2.18}{\xymatrix{ X_{u^*\rho}  \ar[r] \ar[d]_{s_{u^*\rho}} \ar@{}[rd]|\Box& X_\rho \ar[d]^{s_\rho} \\
T\times_k X \ar[r]_{u\times 1 }&  S\times_k X}} 
which makes 1), 2), 3) functorial.

\end{itemize}
\end{thm}
\begin{proof} 
We define $A=(B, D) \subset M=(V, \nabla)$ to be the objects in $\text{\sf Ind-FConn}(X)$ 
which are mapped by $\rho$  to ${p_2}_*\sO_{S\times_k S}\subset s_*\sO_\Pi$ defined in  Theorem \ref{thm2.6}. 
By construction, $B=p_{2*}\sO_{S\times_kX}$, where $p_2:S\times_kX\to X$.
Since ${p_2}_*\sO_{S\times_k S}$ is an $\sO_S$-algebra with trivial action of $\Pi$, 
it follows that $A$, as an object in $\text{\sf Ind-FConn}(X)$, is trivial and
 is an algebra over $(\sO_X, d)$.
 The algebra structure forces via the Leibniz formula  
the connection $D$ to be of the form $p_{2*}(D')$, where $D'$ is a relative connection
on $S\times_kX/S$. The inclusion $(\sO_X,d)\subset (B,D)$ implies that $D'(1)=d(1)=0$, hence
 $D'=d_{S\times_kX/S}$.
Thus we have 
\ga{2.19}{{\rm Spec}_X B=S\times_k X, \ B=(p_2)_*\sO_{S\times_k X}, 
\\ D=(p_2)_*(d_{S\times_kX/S}:\sO_{S\times_kX}\longrightarrow \Omega^1_{S\times_kX/S}) \notag \\ 
\rho((p_2)_*(\sO_{S\times_k X},d_{S\times_kX/S}))= (p_2)_*\sO_{S\times_k S}.\notag
}
We now apply to $M\supset A$ a similar argument. Since $ s_*\sO_\Pi$ is a
${p_2}_*\sO_{S\times_k S} $-algebra as a representation of $\Pi$,
$M$ is an $A$-algebra object in $\text{\sf Ind-FConn}(X)$. 
The geometric information yields 
that $V$ is a $B$-algebra. Thus, setting $X_\rho=\Spec_XV$, we obtain a morphism
$s_\rho:X_\rho\to S\times_kX$. The algebra structure in $\text{\sf Ind-FConn}(X)$ forces
 $\nabla$ to be coming from a relative connection $\nabla':\sO_{X_\rho}\to 
\Omega^1_{X_\rho/S}$.
 Moreover the inclusion $A\subset M$ implies $\nabla'(1)=D(1)=0$, thus $\nabla'=d_{X_\rho/S}$.
 Summarizing, 
Theorem \ref{thm2.6}, 3) translates as follows
\ga{2.20}{\xymatrix{X_\rho={\rm Spec}_X(V)  \ar[dr]_{\pi_\rho}  \ar[r]^{\qquad s_\rho} &  S \times_k X \ar[d]^{p_2}\\
 & X
}\\
V=(\pi_\rho)_*\sO_{X_\rho}, \ \nabla=(\pi_\rho)_* (d_{X_\rho/S}:
\sO_{X_\rho}\longrightarrow \Omega^1_{X_\rho/S}), \notag \\  
\rho((\pi_\rho)_*(\sO_{X_\rho},d_{X_\rho/S}))=s_*\sO_\Pi.\notag
}
Now  Theorem \ref{thm2.6}, 1) together with \eqref{2.6}  shows 1). To show 2) we first notice that 
$\text{\sf Ind-FConn}(X)$ is a full sub-category of the category of flat connections on $X$.
Therefore, by  Theorem \ref{thm2.6} 4), $(B,D)$ is the maximal trivial sub-connection of $(V,\nabla)$ in the latter
category. So one obtains that $H^0_{DR}(X, M)=B$, or said differently,  
$R^0( p_\rho)_{* DR}(X_\rho/S):= R^0(p_\rho)_*(\Omega^\bullet_{X_\rho/S})=\sO_S$.\\[.2cm]
We prove 3).  For $N $ in ${\sf FConn}(X)$, 
one has $s^*\rho(N)\cong t^*\rho(N)$ by \eqref{2.16}. Since $\rho(N)$ is locally free (\cite[1.6]{DeGroth}), projection formula applied to $s$ yields 
$\rho(N)\otimes_{\sO_S} s_*\sO_\Pi\cong s_*t^*\rho(N)$ as representations of $\Pi$. The Tannaka dual of the left  hand side is 
$(\pi_{\rho *}\pi_\rho^*N$ while the Tannaka dual of the right hand side is 
$(\pi_\rho)_* (p_\rho^{-1}\rho(N)\otimes_{p_\rho^{-1}\sO_S} \sO_{X_\rho}, 
1\otimes d_{X_\rho/S})$. This shows 3). \\[.2cm]
For  Claim 4), we just notice that by projection formula applied to the locally 
free bundle $\rho(N)$ and the connection $\pi_\rho^*N$ relative to $S$, 
one has 
\ml{}{
R^0(p_\rho)_{*DR} ((p_\rho^{-1}\rho(N)\otimes_{p_\rho^{-1}\sO_S}
 \sO_{X_\rho}, 1\otimes d_{X_\rho/S}))=\notag \\   
 \rho(N) \otimes_{\sO_S} R^0(p_\rho)_{*DR} (\sO_{X_\rho}, d_{X_\rho/S})
 \notag .}
  By 2), this expression is equal to $\rho(N)$, while by 3), it is isomorphic to
   $R^0(p_\rho)_{*DR}(X_\rho/S, \pi_\rho^*N)$. 
 This shows 4). 
Finally, the  functoriality in 5) is the translation of the base change property \cite[(3.5.1)]{DeGroth}.
\end{proof}
\begin{defn} \label{defn2.8}
\begin{itemize}
\item[1)]
We fix an embedding $k\subset \bar{k}$ which defines ${\rm Spec}(\bar{k})$ as a 
$k$-scheme. Let $\rho: {\sf FConn}(X)\to {\sf Vec}_{\bar{k}}$ 
be a fiber functor. Then $\pi_\rho: X_\rho \to X$ in \eqref{2.13} has the factorization 
\ga{2.21}{\xymatrix{X_\rho  \ar[dr]_{\pi_\rho}  \ar[r]^{s_\rho\qquad} &  
\Spec(\bar{k}) \times_k X \ar[d]^{p_2}\\
 & X
}
}
where $s_\rho$  is a principal bundle over $X$ under the $\bar{k}$-pro-finite group scheme $\Pi^\Delta$. In particular it is a pro-finite \'etale covering. Thus $\pi_\rho$ is an (infinite) \'etale covering of $X$ 
which we call {\it the universal covering of} $X$ {\it associated to} $\rho$. 
\item[2)]Let $\rho=\tau$ be the tautological functor which assigns to a connection $(V,\nabla)$ its 
underlying bundle $V$. Then by definition
\ga{}{\big(X_\tau\xrightarrow{s_\tau} X\times_sX\big)= \big(\Pi \xrightarrow{(t,s)} X\times_kX\big),\notag}
with $\Pi={\rm Aut}^{\otimes }(\tau)$. Since this groupoid plays a special r\^ole, we denote it by $\Pi(X,\tau)$ and call it {\it the 
total fundamental groupoid scheme of} $X$. 

\end{itemize}
\end{defn}
\begin{prop} \label{prop2.9} Let $\rho:{\sf FConn}(X)\longrightarrow {\sf Vec}_k$
 be a neutral fiber functor.  Then for all finite field extensions $L\supset k$, $\rho$ 
lifts uniquely to $\rho_L: {\sf FConn}(X\times_k L)\to {\sf Vec}_L$. 
\end{prop}
\begin{proof} Since $\rho$ is a neutral functor, 
$\Pi={\rm Aut}^{\otimes}(\rho)$ is a group scheme over $k$. Tannaka duality reads 
${\sf FConn}(X)\xrightarrow{\rho} {\rm Rep}(\Pi)$.  
By Property \ref{eig2.4}, 4), ${\sf FConn}(X\times_k L)$ is the $L$-base change of
${\sf FConn}(X)$ in the sense of \cite[Corollary~5.6]{Hai}.
So Tannaka duality lifts to the duality $\rho_L:{\sf FConn}(X)\to \Rep_L(\Pi\times_k L)$
which fits in the following commutative diagram
\ga{2.22}{\xymatrix{{\sf FConn}(X)\ar[rr]^\rho_\cong\ar[d]& & \Rep_k(\Pi)\ar[d] \\
{\sf FConn}(X\times_k L)\ar[rr]_{\rho_L}&& \Rep_L(\Pi\times_k  L).
}}
We define the functor $\rho_L$ as follows.
Let $\pi$ be the projection $X \times_k L \to X$. Then for any connection $N\in{\sf FConn}(X\times_k L)$
one has the canonical projection $\pi^*\pi_*N\to N$ of connections in ${\sf FConn}(X \times_k L)$.
Consequently, by considering the dual connection $N^\vee$, one obtains an injection
$N\to (\pi^*\pi_*N^\vee)^\vee\cong \pi^*(\pi_*N^\vee)^\vee$. Thus $N$ is uniquely determined
by a morphism $f_N\in H^0_{DR}(X_L,\pi^*((\pi_*N)^\vee\otimes (\pi_*N^\vee)^\vee ))$.
Via the base change isomorphism
\ga{2.23}{H^0_{DR}(X\times_k L,\pi^*M)\stackrel\cong\longrightarrow H^0_{DR}(X,M)\otimes_kL,}
one can interpret $f_N$ as an element of 
$H^0_{DR}(X,(\pi_*N)^\vee\otimes (\pi_*N^\vee)^\vee )\otimes_kL\cong \Hom_{{\sf FConn}(X)}
\pi_*N,(\pi_*N^\vee)^\vee)\otimes_k L$.\\[.2cm]
We define now the functor $\rho_L$ on objects of ${\sf FConn}(X \times_k L)$ of the form
$\pi^*M$ by setting $\rho_L(\pi^*M)=\rho(M)\otimes_kL$ and  
 \eqref{2.23} allows us to define $\rho_L$ on morphisms from $\pi^*M$ to $\pi^*M'.$ \\[.2cm]
Next, we extend $\rho_L$ to an arbitrary object $N$ by defining $\rho_L(N)$ to be
the image of $\rho_L(f_N)$ in $\rho_L((\pi_*N^\vee)^\vee)$. \\[.2cm]
As for unicity, if $\rho'$ is another lifting of $\rho$, then it has to agree on $f_N$ with $\rho_L$, thus has to agree with $\rho_L$ for all $N$. 
\end{proof}
\begin{nota} \label{nota2.10}
With notations as in Proposition \ref{prop2.9}, we denote by $\rho_L$ the unique lifting of $\rho$ to $X\times_k L$ and by $\rho\times_k\bar{k}$
the induced lifting on $X\times_k \bar{k}$, once an algebraic closure $k\subset \bar{k}$ of $k$ has been fixed.

\end{nota}
\begin{thm} \label{thm2.11}
Let $X$ be a smooth scheme of finite type defined over a field $k$ of characteristic 0. 
We fix an embedding $k\subset \bar{k}$ which defines ${\rm Spec}(\bar{k})$ as a $k$-scheme. 
Assume there is a neutral fiber functor  $\rho: {\sf FConn}(X)\to {\sf Vec}_{k}$. Then one has an isomorphism
\ga{2.24}{\xymatrix{ \ar[dr]_{s_\rho\times_k \bar{k}}  X_\rho \times_k\bar{k} \ar[rr]^{\cong} & &  X_{\rho\times_k \bar{k}} \ar[dl]^{s_{\rho\times_k \bar{k}}}\\
& X\times_k \bar{k}=\bar{k}\times_k X 
}}
In particular, this yields a 
cartesian diagram
\ga{2.25}{
\xymatrix{ X_{\rho\times_k \bar{k}}\cong \bar{k}\times_k X_\rho 
 \ar[d] \ar[r]
\ar@{}[rd]|\Box
& \bar{k}\times_k X \ar[d]^{p_2} \\
X_\rho  \ar[r]_{s_\rho} &  X }
.} 
\end{thm}
\begin{proof}
If $\Pi={\rm Aut}^{\otimes }(\rho)$ is as in the proof of Proposition \ref{prop2.9}, then by construction, $\Pi\times_k \bar{k}={\rm Aut}^{\otimes }(\rho)\times_k \bar{k}$. Thus base change as in Theorem \ref{thm2.7}, 5) implies the wished base change property for $X_\rho$.

\end{proof}
\noindent We now come back to the converse of Properties \ref{eig2.4}, 3). 
\begin{cor} \label{cor2.12} Let $X$ be as in Definition \ref{defn2.1}. Then
 ${\sf FConn}(X)$ is a semi-simple category. 
Consequently any object of ${\sf FConn}(X)$ is Weil-finite.
\end{cor}
\begin{proof}
We have to show that an exact sequence $\epsilon:0\to M\to P\to N \to 0$ splits,
or equivalently, that the corresponding cohomology class in $H^1_{DR}(X,N^\vee\otimes M)$, 
also denoted by $\epsilon$,
 vanishes. Since de Rham cohomology fulfills
 base change, $\epsilon$ vanishes if and only if $\epsilon\otimes_k\bar k$
vanishes in $H^1_{DR}(X\times_k\bar k,N^\vee\otimes M)$. Thus we may assume
that $k=\bar k$. Then we can find a $k$-rational point $x\in X(k)$ and ${\sf FConn}(X)$
is a (neutral) Tannaka category with respect to the fiber functor $\rho_x((V,\nabla))=V|_x$. 
The exact sequence $\epsilon: 0\to M\to P\to  N \to 0$ lies in the category $\langle M\oplus N\oplus P\rangle$,
which is finite by Properties \ref{eig2.4}, 3).  Hence its Tannaka group $H$ is a finite
 group scheme over $\bar{k}$, and one has a surjective
factorization 
\ga{2.26}{\xymatrix{ X_\rho \ar@{->>}[r] \ar[dr]_{\pi_\rho} & X_H \ar[d]^{\pi_H}\\
 & X}
}
where $\pi_H$ is a principal bundle under $H$. As $\pi_H^*$ is injective on de Rham cohomology, we are reduced to the case where $N=M=(\sO_X,d), P\cong(\sO_X,d) \oplus (\sO_X,d)$. Then necessarily $\epsilon$ splits.

\end{proof}
\begin{rmks} \label{rmk2.13}\begin{itemize}
\item[1)] The notations are as in Theorem \ref{thm2.7}. Assume that
  $\rho$ is a neutral fiber functor. 
Then  $\Pi={\rm Aut}^{\otimes }(\rho)$ is an affine pro-finite $k$-group scheme, $S={\rm Spec}(k)$  and by definition $s_*\sO_\Pi=k[\Pi]$. 
Then \eqref{2.20} yields in particular 
\ga{2.27}{\rho((\pi_\rho)_*(\sO_{X_\rho}, d_{X_\rho/k}))=k[\Pi].}
Let us assume that $\rho=\rho_x$, where $x\in X(k)$ and $\rho_x((V,\nabla))=V|_x$. Then \eqref{2.27} reads
\ga{2.28}{(\pi_{\rho_x})_*(\sO_{X_{\rho_x}})|_x=k[\Pi].}
Here the $k$-structure on the $k$-group scheme $\Pi$ is via the residue field $k$ of $x$. Thus in particular, 
\ga{2.29}{1\in \Pi(k)=\pi_{\rho_x}^{-1}(x)(k)\subset X_{\rho_x}(k)}
and the rational point $x$ of $X$ lifts all the way up to the universal covering of $X$ associated to $\rho_x$.  \\[.2cm]
\item[2)] Theorem \ref{thm2.7}, 4) says in loose terms that all fiber functors of ${\sf FConn}(X)$ are cohomological, as they are canonically isomorphic to a 0-th relative de Rham cohomology. On the other hand, there are not geometric in the sense that they do not come from rational points of $X$ or of some compactification (for the latter, see discussion in section 5). A simple example is provided by a smooth projective rational curve 
 defined over a number field $k$, and without any rational point. Then ${\sf FConn}(X)$ is trivial, thus has the neutral fiber functor $M=(V,\nabla)\mapsto H^0_{DR}(X, V)$. Yet there are no rational points anywhere around. 
\end{itemize}
\end{rmks}
\noindent Recall the definition of Nori's category of Nori finite bundles (\cite[Chapter~II]{N2}). 
Let $X$ be a proper reduced scheme defined over
a perfect field $k$. Assume $X$ is connected in the sense that $H^0(X, \sO_X)=k$. 
Then the category $\sC^N(X)$ of Nori finite bundles is the full sub-category of the 
quasi-coherent category ${\sf QCoh}(X)$ consisting of bundles which in ${\sf QCoh}(X)$ 
are sub-quotients of Weil-finite bundles, where a bundle is Weil-finite when there are 
polynomials $f,g \in \N[T], f\neq g$ such that $f(M)$ is isomorphic to $g(M)$. 
\\[.2cm]
Let $\omega:\sC^N(X)\longrightarrow
{\sf QCoh}(X)$ be a fiber functor, where $S$ is a scheme over $k$. 
Then one can redo word by word  the whole construction of Theorem \ref{thm2.7} with 
$({\sf FConn}(X), \rho)$ replaced by  $(\sC^N(X), \omega)$, as it goes purely 
via Tannaka duality. We denote by $\Pi^N={\rm Aut}^{\otimes}(\omega)$ 
the $k$-groupoid scheme acting transitively on $S$. One obtains the following. 
\begin{prop} \label{prop2.14} Let $X$ be a connected proper reduced scheme of finite type 
over a perfect field $k$. Let $\omega$ be a fiber functor $\sC^N(X)\longrightarrow
{\sf QCoh}(S)$ where $S$ is a $k$-scheme and let  $\Pi^N={\rm Aut}^{\otimes}(\omega)$ be the
corresponding Tannaka groupoid scheme acting on $S$. 
  Then there is a diagram of  $k$-schemes
  \ga{2.30}{\xymatrix{& \ar[ddl]_{p^N_\omega} X^N_\omega \ar[d]^{s^N_\omega} \ar[ddr]^{\pi^N_\omega}\\
& S\times_k X \ar[dl]^{p_1}  \ar[dr]_{p_2}\\
S &  & X
}
}
with the following properties.
\begin{itemize}
\item[1')]$s^N_\omega$ is a 
$(\Pi^N)^\Delta$-principal bundle, that is
$$X^N_\omega \times_{S\times_k X}X^N_\omega\cong(\Pi^N)^\Delta \times_S X^N_\omega .$$ 
\item[2')] $R^0 (p^N_\omega)_*(\sO_{X^N_\omega})=\sO_S$. 
\item[3')] For all objects $V$ in $\sC^N(X)$, the bundle
$(\pi^N_\omega)^*V$  is endowed with an isomorphism with the bundle 
$(p^N_\omega)^*\omega(V)$ which is compatible with all morphisms in 
$\sC^N(X)$. 
\item[4')]
 One recovers the fiber functor $\omega$ via $X^N_\rho$ by an isomorphism
$\omega(M)\cong R^0(p_\omega^N)_{*}(X_\omega/S, \pi_\omega^*M)$ which  is compatible with all morphisms in $\sC^N(X)$. In particular, the data in \eqref{2.30} are equivalent to the datum $\omega$ (which defines $\Pi^N$). 
\item[5')] This construction is compatible with base change: if $u: T\to 
S$ is a morphism of $k$-schemes, and $u^*\omega: \sC^N(X)\to {\sf QCoh}(T)$ is the composite fiber functor, then one has a cartesian diagram
\ga{2.31}{\xymatrix{ X^N_{u^*\omega}  \ar[r] \ar[d]_{s^N_{u^*\omega}} \ar@{}[rd]|\Box& X^N_\omega \ar[d]^{s^N_\omega} \\
T\times_k X \ar[r]_{u\times 1 }&  S\times_k X}} 
which make 1), 2), 3) functorial.
\item[6')] (see Definition \ref{defn2.8}, 2). If $\omega$ is the tautological fiber functor $\iota$ defined by  $\iota(V)=V$, then 
\ga{}{\big(X_\iota^N\xrightarrow{s_\tau} X\times_sX\big)= \big(\Pi^N \xrightarrow{(t,s)} X\times_kX\big).\notag}

\end{itemize}
 \end{prop}

\begin{thm} \label{thm2.15}
Let $X$ be a smooth proper scheme of finite type over a field $k$ of characteristic 0 
with $k=H^0(X, \sO_X)$. 
Then the functor $F:{\sf FConn}(X)\to \sC^N(X)$, $F((V,\nabla))=V$ 
 is an equivalence of Tannaka categories.
\end{thm}
\begin{proof} 
Since both categories ${\sf FConn}(X)$ and $\sC^N(X)$ satisfy the base change
property in the sense of Property \ref{eig2.4}, 4), we may assume that
$k=\bar k$.\\[.2cm]
Let $V$ be a finite bundle, we want to associate to it a connection $\nabla_V$,
such that $(V,\nabla_V)$ is a finite connection. Denote by $\sT$ the  full tensor
sub-category of $\sC^N(X)$ generated by $V$. It is a finite category.
We apply the construction in Proposition \ref{prop2.14} to $\omega=\iota:
\sT\longrightarrow {\sf QCoh}(X)$, defined by $\iota(V)=V$. Then \eqref{2.30} reads
\ga{2.32}{
\xymatrix{& \ar[ddl]_{t} X_\sT=\Pi_\sT \ar[d] \ar[ddr]^{s}\\
& X\times_k X \ar[dl]^{p_1}  \ar[dr]_{p_2}\\
X &  & X
}}
where $t,s$ are the structure morphisms. The scheme 
$X_\sT$ is a principal bundle over $X\times_kX$ under the finite $X$-group scheme
 $(\Pi_\sT)^\Delta$.
For an object $V$ in $\sT$, 3') becomes a functorial isomorphism
\ga{2.33}{t^*V\cong s^*V.}
Let $d_\text{rel}$ denote the relative differential $\sO_{X_\sT}\to \Omega^1_{X_\sT}/s^*\Omega^1_X$  on $X_\sT/X$ with respect to the
morphism $s$. Then the bundle $s^*V$ carries the canonical connection
\ga{2.34}{s^*V=s^{-1}V\otimes_{s^{-1}\sO_X} \sO_{X_\sT}
\xrightarrow{\id_V\otimes d_\text{rel}}
s^{-1}V\otimes_{s^{-1}\sO_X} (\Omega^1_{X_\sT}/s^*\Omega^1_X).}
Then \eqref{2.33} implies that \eqref{2.34} can be rewritten as
\ga{2.35}{
t^*V \xrightarrow{\id_V\otimes d_\text{rel}} t^*(V)\otimes_{\sO_{X_\sT}}t^*(\Omega^1_X).}
Applying $R^0t_*$ to \eqref{2.35}, the property 2') together with projection formula 
implies that one obtains a connection
\ga{2.36}{ \nabla_V:=R^0t_*(\id_V\otimes d_\text{rel}): V\to V \otimes_{\sO_X} \Omega^1_X.}
Integrability of $\id_V\otimes d_\text{rel}$ implies integrability of $\nabla_V$. 
The compatibility of $\id_V\otimes d_\text{rel}$ with morphisms in $\sC^N(X)$ implies the
compatibility of $\nabla_V$ with morphisms in $\sC^N(X)$ as well. We conclude 
that $\nabla_V$ defines a functor
\ga{2.37}{\sC^N(X)\to {\sf Conn}(X), \ V\mapsto (V, \nabla_V).}
It remains to show that $(V,\nabla_V)$ is in ${\sf FConn}(X)$. To this aim, we
fix a ${k}$-rational point $x:\Spec(k) \to X$ and consider the restriction of \eqref{2.32}
to $X_{\sT,x}\to X\times_k x$, where $X_{\sT,x}=s^{-1}x$. The restriction of 
$t$ to $X_{\sT,x} $ 
will be denoted by $t_x: X_{\sT,x}\to X$. This is a principal bundle under the finite 
$k$-group scheme $(\Pi_\tau)^\Delta|_x$, the fiber of $(\Pi_\tau)^\Delta\to X$ above $x$. 
Since the connection $\id_V\otimes d_\text{rel}$ in \eqref{2.34}, or equivalently \eqref{2.35}
 is relative to the $X$-factor on the right, we can restrict it to $X_{\sT,x}$ to obtain 
 a connection 
$$(\id_V\otimes d_{{\rm rel}})|_{X_{\sT,x}}: t_x^*V \xrightarrow{\id_V\otimes d_\text{rel}}
  t_x^*(V)\otimes_{\sO_{X_{\sT,x}}}t_x^*(\Omega^1_X).$$
Since $t^*V$ was generated by the relatively flat sections $s^{-1}V$ of $\id_V\otimes d_{{\rm rel}}$, 
$ t_x^*V$ is generated by the relative flat sections 
$V|_x$ of  $(\id_V\otimes d_{{\rm rel}})|_{X_{\sT,x}}$. 
On the other hand, by construction
\ga{2.38}{t_x{}^*(V,\nabla_V)=
(\id_V\otimes d_\text{rel})|_{X_{\sT,x}}.}
Thus $(V,\nabla_V)$ is trivializable when pulled back to the principal bundle $t_x: X_{\sT,x}\to X$ 
under the finite group scheme $(\Pi_\tau)^\Delta|_x$. 
Hence $(V,\nabla_V)$  is a sub-connection of a direct sum of copies of
 the connection $(W,\nabla_W):=
(t_x)_*(\sO_{X_{\sT,x}}, d_\text{rel})$, which is known
to be finite as it fulfills the equation 
\ga{2.39}{ (W,\nabla_W)^{\otimes 2}\cong (W,\nabla_W)^{\oplus {\rm rank} W}.}
Hence $(V,\nabla)$ is finite.
We conclude that \eqref{2.37} has image in ${\sf FConn}(X)$. Since the composite functor
\ga{2.40}{{\sf FConn}(X)\xrightarrow{F\ \inj} \sC^N(X) \xrightarrow{\eqref{2.37}\ \inj} {\sf FConn}(X)}
is the identity, and both $F$ and \eqref{2.37} are full. This finishes the proof.
\end{proof}

\section{Splitting of  groupoid schemes}
\begin{notas} \label{nota3.1} Throughout this section, we will use the following notations. 
$k$ is a field of characteristic 0, endowed with an algebraic closure $\bar{k}\supset k$. 
$\sC$ is an abelian $k$-linear rigid tensor category, endowed with a fiber functor $\rho: \sC\to {\sf Vec}_{\bar{k}}$. 
By \cite[Th\'eor\`eme 1.12]{DeGroth}, the groupoid scheme $\Pi={\rm Aut}^{\otimes}(\rho)$ 
over $k$, acting on ${\rm Spec}(\bar{k})$, via $(t,s): \Pi\to 
{\rm Spec}(\bar{k})\times_k {\rm Spec}(\bar{k})$, is affine and acts transitively.
 $\sC$ is equivalent via the fiber functor to the category of finite dimensional
representations of $\Pi$:
\ga{3.1}{ \sC\xrightarrow{\rho \ \cong} {\rm Rep}_f(\Spec(\bar{k}): \Pi).}
As in \eqref{2.3}, we denote by $\Pi_s$ the scheme $\Pi$ viewed as a $\bar{k}$-scheme via 
the projection $s: \Pi\to {\rm Spec}(\bar{k})$. 
So 
\ga{3.2}{\Pi_s(\bar{k})=\{\xymatrix{ {\rm Spec}(\bar{k})\ar[r] \ar[dr]_{{\rm Id}} & \Pi \ar[d]^s\\
& {\rm Spec}(\bar{k})\}.}
}
For $\sC$ the trivial category, the objects of which being sums of finitely many copies of the unit object, 
$\Pi$ is the groupoid scheme  ${\rm Spec}(\bar{k})\times_k {\rm Spec}(\bar{k})$ 
while $\Pi_s(\bar k)$ is  identified via Galois theory with the pro-finite group  ${\rm Gal}(\bar{k}/k)$.   
\end{notas}

\begin{thm} \label{thm3.2}
\begin{itemize}
\item[1)] There exists a group structure on $\Pi_s(\bar k)$
such that the map 
$$(t,s)|_{\Pi_s}: \Pi_s(\bar k)\to ({\rm Spec}(\bar{k})\times_k {\rm Spec}(\bar{k}))_s(\bar k)\cong 
{\rm Gal}(\bar{k}/k)$$
 is a group homomorphism. 
\item[2)]Splittings of $(t,s)|_{\Pi_s}: \Pi_s(\bar k)\to ({\rm Spec}(\bar{k})\times_k 
{\rm Spec}(\bar{k}))_s(\bar k)\cong {\rm Gal}(\bar{k}/k)$ 
as group homomorphisms are in one to one correspondence  with splittings of 
$(t,s): \Pi\to {\rm Spec}(\bar{k})\times_k {\rm Spec}(\bar{k})$ as $k$-affine groupoid 
scheme homomorphisms. 
\item[3)]  The projection $(t,s): \Pi\to {\rm Spec}(\bar{k})\times_k {\rm Spec}(\bar{k})$ has a 
section of  groupoid schemes over $k$ if and only if $\sC$ has a neutral fiber functor.
 More precisely,
there is a 1-1 correspondence between neutral fiber functors of $\sC$ and splittings of $(t,s)$
up to an inner conjugation of $\Pi$ given by an element of $\Pi^\Delta(\bar k)$.
\end{itemize}
\end{thm}
\begin{proof}
 Let $\theta \in \Pi_s(\bar{k})$, set
 $\gamma=t\theta:{\rm Spec}(\bar{k})\stackrel\theta\to \Pi_s\stackrel t\to\Spec(\bar k)$.
 As the morphisms are all over ${\rm Spec}(k)$, $ \gamma$ is an element of $ {\rm Gal}(\bar{k}/k)$. 
 By definition, $\theta$ is a tensor  isomorphism between $\rho$ and 
 $\gamma^*(\rho)$, where $\gamma^*(\rho)$ is the ${\rm Spec}(\bar{k})$-valued fiber functor 
 on $\sC$ which is the composite of $\rho: \sC \to {\sf Vec}_{\bar{k}}$ with 
 $\gamma^*: {\sf Vec}_{\bar{k}}\to {\sf Vec}_{\bar{k}}, \ V\mapsto V\otimes_{\gamma} \bar{k}$:
\ga{3.3}{ \rho \xrightarrow{\theta \ \cong}  \gamma^*\circ \rho.}
Thus given $\theta'$ with image $\gamma'$, the group structure is simply
given by the compositum 
\ga{3.4}{\theta \cdot \theta': \rho \xrightarrow{\theta' \ \cong}  (\gamma')^*\circ \rho
\xrightarrow{ (\gamma')^* \circ \theta  \cong } (\gamma')^*\circ \gamma^* \circ \rho= 
(\gamma\circ \gamma')^* \circ \rho. }
The inverse to $\theta$ is $\gamma^*(\theta^{-1})$.
This shows 1). \\[.2cm]
Assume that one has a splitting $\sigma_s:{\rm Gal}(\bar k/k)\to \Pi_s(\bar k)$
of the map $(t,s)|_{\Pi_s(\bar{k})}:\Pi_s(\bar{k})\to {\rm Gal}(\bar k/k)$. 
 This means that to  any $\gamma$
in ${\rm Gal}(\bar k/k)$, one assigns a map $\sigma_s(\gamma):\Spec(\bar k)\to \Pi_s(\bar{k})$, such that
$$ \sigma_s(\gamma)\cdot\sigma_s(\gamma')=\sigma_s(\gamma\circ \gamma'); \quad t(\sigma(\gamma))=\gamma.$$
As in the discussion above, $\sigma_s(\gamma)$ defines a natural isomorphism
$\sigma_s(\gamma):\rho\to \rho\otimes_\gamma\bar k$. For a $\bar k$-vector space $V$ 
we identify $V\otimes_\gamma\bar k$
with $V$ by setting $v\otimes_\gamma a\mapsto v\gamma(a)$. Consequently $\sigma_s(\gamma)$
yields a natural action of ${\rm Gal}(\bar k/k)$ on $\rho$:
\ga{3.5}{\gamma:\rho \stackrel{\sigma_s(\gamma)}\longrightarrow \rho 
\otimes_\gamma\bar k\cong\rho,
\ \gamma\in{\rm Gal}(\bar k/k). }
Here ``natural action'' means that it commutes with morphism in $\sC$. Galois theory 
applied to the values of the fiber functor $\rho$ 
implies that there exists a $k$-form $\rho_0$ for $\rho$: $\rho=\rho_0\otimes_k\bar k$.
Now the Tannaka construction of $\Pi$ from $\rho$ \cite[Section~4.7]{DeGroth} tells us that there exists
a homomorphism of groupoid schemes acting on $\Spec(\bar k)$: $\Pi\longrightarrow 
\Spec(\bar k)\times_k\Spec(\bar k)$.
The converse claim is trivial. This finishes the proof of   2).\\[.2cm]
We turn now to the proof of 3). Notice that the structure map $(t,s)$ considered as a homomorphism
of groupoid schemes corresponds through Tannaka duality to
 the tautological fully faithful functor $\beta:{\sf Vec}_{k,f}\to \sC$,
$V\mapsto V\otimes_kI$, which has the property that $\rho\circ \beta$ is the base change $\tau: {\sf Vec}_{k,f}\to {\sf Vec}_{\bar{k}}, \ V\mapsto V\otimes_k \bar{k}$. Here ${\sf Vec}_{k,f}$ means the category of finite
dimensional $k$-vector spaces.  (See \cite[Proof of Theorem~2.11]{DeMil}). 
Then splittings of $(t,s)$ as homomorphisms of $k$-groupoid schemes acting on ${\rm Spec}(\bar k)$
are in one to one  correspondence with splittings $\sigma$ of the  functor $\beta$ which are compatible
with $\tau$ and $\rho$. This means one has 
a functor $\sigma: \sC \to {\sf Vec}_{k,f}$ such that $ \sigma\circ \beta={\rm Id}$, together with
an isomorphism of tensor functor $d: \tau\circ \sigma\to \rho$. So the following diagram 
\ga{}{\notag\xymatrix@W=5pt{
{\sf Vec}_{k,f} \ar[rd]_\tau &&\sC\ar[dl]^\rho \ar[ll]_\sigma\\ &{\sf Vec}_{\bar k}&
}}
commutes up to $d$.  
\\[.2cm]
By \cite[Proposition~8.11]{DeGroth}, the isomorphism functor ${\rm Iso}^{\otimes}(\rho,\tau\circ \sigma)$
is representable by a torsor over ${\rm Spec}(\bar k)$ under the $\bar k$-group scheme $\Pi^\Delta$. 
Thus, through Tannaka duality, neutral fiber functors of $\sC$ are in one to one correspondence
with splittings to $(t,s)$ up to an inner conjugation of $\Pi$ given by an element of 
$\Pi^\Delta(\bar k)$. This finishes the proof of 3).

\end{proof}

\section{Grothendieck's arithmetic fundamental group}
For a scheme $X$ define over a field $k$ Grothendieck introduces in \cite[Section~5]{Groth}   the category 
${\sf ECov}(X)$ of finite \'etale coverings $\pi:Y\to X$, with Hom-Sets being $X$-morphisms.
The choice of a geometric point $\bar x$ of $X$ defines a fiber functor $\pi\mapsto \pi^{-1}(\bar x)$
with value in the category of finite sets. The automorphism group of this functor
is called the {\it arithmetic fundamental group of $X$ with base point  $\bar x$} and is denoted by
$\pi_1(X,\bar x)$. It is an abstract
group, endowed  with the pro-finite topology stemming from all its finite quotients.
The main theorem claims an equivalence between finite sets with continuous $\pi_1(X,\bar x)$-action  and finite \'etale coverings of $X$. The equivalence also extends to an
equivalence between pro-finite sets with continuous  $\pi_1(X,\bar x)$-action 
and pro-finite \'etale coverings of $X$. In particular, the action of $\pi_1(X,\bar x)$ on
itself via translations defines the {\em universal pro-finite \'etale covering of $X$ based at $\bar x$.} It will be denoted by
$ \tilde{\pi}_{\bar x}:  \tilde X_{\bar x}\to X$.
By definition, $\tilde{\pi}_{\bar x}^{-1}(\bar x)=\pi_1(X,\bar x)$. \\[.2cm]
This section is devoted to the comparison between Grothendieck's arithmetic fundamental
group and our fundamental groupoid scheme, as well as between the universal pro-finite
\'etale covering and a  special case of the covering constructed in Theorem \ref{thm2.7}. 
\begin{notas}\label{nota4.1} In this section, $X$ is again a smooth scheme of finite type over a field 
$k$ of characteristic 0, with the property $k=H^0_{DR}(X)$.
For a fiber functor $\rho$ of ${\sf FConn}(X)$ we denote by $$\Pi(X,\rho)={\rm Aut}^{\otimes}(\rho)$$ the corresponding
Tannaka $k$-groupoid scheme. (As compared to the notations of section 2 where we simply used the notation $\Pi$, we emphasize here $\rho$). For a geometric point $\bar x:\Spec(\bar k)\to X$, we denote by  $\rho_{\bar x}
:{\sf FConn}(X)\to {\sf Vec}_{\bar k}$ 
the fiber functor that assigns to a connection the fiber of the underlying bundle at $\bar x$. We simplify the notation by setting 
 $$\Pi(X,\bar x):=\Pi(X,\rho_{\bar x}).$$ We call this $k$-groupoid scheme acting on ${\rm Spec}(\bar{k})$ 
 the {\em fundamental groupoid scheme of $X$ with base point $\bar x$}.\\
Recall that the embedding $k\subset \bar{k}$ here is defined by the residue field of $\bar{x}$. Let $L$ be a finite field extension of $k$ in $\bar k$.
 We define the $L$-base change $\Pi(X,\bar x)_L$
of $\Pi(X,\bar x)$ by the following cartesian product
$$\xymatrix{
\ar@{}[dr] |\square
\Pi(X,\bar x)_L\ar[r]^q \ar[d]&\Pi(X,\bar x)\ar[d]\\
\Spec L\ar[r]_{\Delta\qquad\quad}& \Spec(L)\times_k\Spec(L)
}$$
where the morphism $\Pi(X,\bar x)\to \Spec(L)\times_k\Spec(L)$ is the composition of
$(t,s)$ with the projection $\Spec(\bar k)\times_k\Spec(\bar k)\to \Spec(L)\times_k\Spec(L)$.
Then $\Pi(X,\bar x)_L$ is an $L$-groupoid scheme acting on $\Spec(\bar k)$ and we have the following
commutative diagram
\ga{4.1}{\xymatrix{ \Pi(X, \bar x)_L \ar[r] \ar[d]_{(t,s)} & 
\Pi(X, \bar x)\ar[d]_{(t,s)} \ar[r] \ar[d]_{(t,s)} &
 {\rm Spec}(L)\times_k {\rm Spec}(L) \ar[d]_{{\rm Id}} \\
{\Spec}(\bar{k})\times_L {\rm Spec}(\bar{k}) \ar[r]  & 
{\Spec}(\bar{k})\times_k {\rm Spec}(\bar{k}) \ar[r] & 
{\rm Spec}(L)\times_k {\rm Spec}(L)}
}
which is an exact sequence of groupoid schemes in the sense of \cite{Hai}.\\[.2cm]
\end{notas}
\begin{lem} \label{lem4.2} 
For all finite extensions $L\supset k$ with $  L\subset \bar{k}$, 
one has a canonical isomorphism 
$${\Pi(X\times_k L,\bar x)\cong \Pi(X,\bar x)_L}$$
which implies that  $\Pi(X\times_k \bar k,\bar x)\cong \Pi(X,\bar x)^\Delta$.
 Consequently the following diagram is an exact sequence of groupoid schemes
\ga{4.2}{\xymatrix{ \Pi(X\times_k\bar k, \bar x)\ar[d] \ar[r]& \Pi(X, \bar x)\ar[d]_{(t,s)}
 \ar[r]& {\rm Spec}(\bar k)\times_k {\rm Spec}(\bar k)  \ar[d]_{{\rm Id}}  \\
 {\rm Spec}(\bar{k}) \ar[r]_{{\rm \Delta\qquad}} &  {\Spec}(\bar{k})\times_k {\rm Spec}(\bar{k})
 \ar[r]_{\id} & 
{\rm Spec}(\bar{k})\times_k {\rm Spec}(\bar{k}) ,}}
where $\Delta$ is the diagonal embedding.
\end{lem}

\begin{proof}
We apply the base change property \ref{eig2.4}, 4) and  \cite[Corollary~5.11]{Hai}
to obtain 
$  \Pi(X\times_k L,\bar x)\cong \Pi(X,\bar x)_L$. Taking the limit on all $L$, one obtains 
$\Pi( X\times_k\bar k, \bar x)\cong \Pi(X\times_k\bar k, \bar x)^{\Delta}$. 
\end{proof}
\noindent
In the sequel, we use the simpler notation $\bar X:=X\times_k \bar k$.
\\[.2cm]
Applying the construction of Theorem \ref{thm2.7} to the category
${\sf FConn}(\bar X)$ equipped with the fiber functor at $\bar x$,  ones obtains
an \'etale cover 
\ga{4.3}{\pi_{\rho_{\bar{x}} }: (\bar X)_{\rho_{\bar x}}\longrightarrow  \bar X}
\begin{lem}\label{lem4.3}
The covering $\pi_{\rho_{\bar{x}} }: (\bar X)_{\rho_{\bar x}}\longrightarrow \bar X$
 is the universal pro-finite \'etale covering $\tilde{\pi}_{\bar{x}}: {\widetilde{(\bar X})}_{\bar x} \to \bar X$ based at $\bar x$. 
 In particular, one has the  identification
$\bar{\kappa}:  \Pi(\bar X,\bar x)(\bar k) \xrightarrow{=} \pi_1(\bar X,\bar x)$ which decomposes as 
\ga{4.4}{\bar\kappa:\Pi(\bar X,\bar x)(\bar k)=\pi_{\rho_{\bar{x}} }^{-1}(\bar{x})=
\tilde{\pi}_{\bar{x}}^{-1}(\bar{x})=
\pi_1(\bar X,\bar x).}
\end{lem}
\begin{proof} Since $\bar X$ is defined over an algebraically closed field $\bar k$,
its universal pro-finite \'etale covering based at $\bar x$ is the pro-limit of
finite Galois coverings $Y\stackrel G\longrightarrow \bar X$ with  $H^0_{DR}(Y)=\bar k$,
which over $\bar k$ simply means $Y$ is connected.
On the other hand, according to Theorem \ref{thm2.7}, $\pi_{\rho_{\bar{x}} }$
 is a principal bundle under the $\bar k$-group scheme
$\Pi(\bar X,\bar x)$. In fact, any finite full tensor sub-category $\sT$ of ${\sf FConn}(X)$
defines a principal bundle $(\bar X)_{\sT,\rho_{\bar x}}$ under the $\bar k$-group
scheme $\Pi_\sT(\bar X,\bar x)$, hence a Galois covering
of $X$ under the finite group $\Pi_\sT(\bar X,\bar x)(\bar k)$.
Also according to Theorem \ref{thm2.7},~2), $H^0_{DR}((\bar X)_{\sT,\rho_{\bar x}})=\bar k$, thus 
 $  (\bar X)_{\sT,\rho_{\bar x}} $ is connected.
 Since ${\sf FConn}(\bar X)$ is the union of its
finite sub-categories, $ \pi_{\rho_{\bar{x}} }: (\bar X)_{\rho_{\bar x}}\to \bar X$ is the pro-limit of the
$\pi_{\sT, \rho_{\bar{x}} }: (\bar X)_{\sT,\rho_{\bar x}}\to \bar X$. Furthermore, by construction,
$\pi_{\sT, \rho_{\bar{x}} }^{-1}(\bar x)= \Pi_\sT(\bar X,\bar x)$ as a $\bar{k}$-group scheme, 
where $\bar{k}$ is the residue field of $\bar x$.  
\\[.2cm]
Conversely, let $p:Y\to X$ be a Galois  covering with Galois group $G$, with $\pi_1(X, \bar x)\surj G$,
 thus  with $Y$ connected. Considering $G$ as a $\bar k$-algebraic  constant group, then $p:Y\to X$ is a principal
bundle under $G$ with $G=p^{-1}(\bar x)$. Then  $M:=p_*(\sO_Y,d)$ is finite as it fulfills the relation
$M^{\otimes 2}\cong  M^{\oplus {\rm deg}(p)}$.  
If we denote by $\sT$ the full tensor sub-category generated by $p_*(\sO_Y,d)$
then $G\cong  \Pi_\sT(\bar X,\bar x)(\bar k)$ and $Y\cong (\bar X)_{\sT,\rho_{\bar x}}$.
Thus $\pi_{\rho_{\bar{x}} }$ is the universal pro-finite \'etale covering of $\bar X$ based at $\bar x$. 
This shows \eqref{4.4} and finishes the proof.
\end{proof}

\begin{thm} \label{thm4.4}
Let $X/k$ be smooth scheme with $H^0_{DR}(X)=k$. 
Let $\bar x\to X$ be a geometric point with residue field $\bar{k}$.  
Then \eqref{4.2} induces an exact sequence of pro-finite groups
\ga{4.5}{1\to \Pi(\bar X,{\bar x})(\bar{k})\to \Pi(X,{\bar x})_s(\bar k) 
\to ({\rm Spec}(\bar k)\times_k {\rm Spec}(\bar k))_s(\bar{k}) 
\to 1 .}
Furthermore, the identity $\bar{\kappa}$ of Lemma \ref{lem4.3} 
extends to an identity of exact sequences of pro-finite groups
\ga{4.6}{\xymatrix{ 1\ar[r] & \Pi(\bar X,\bar x)(\bar k)\ar[r] \ar[d]_{\bar{\kappa}}^{=} & 
 \Pi(X,\bar x)_s(\bar k) \ar[d]_{\kappa }^{=}\ar[r] &
({\rm Spec}(\bar k)\times_k {\rm Spec}(\bar k))_s(\bar{k}) 
\ar[r] \ar[d]^{= }& 1\\
1\ar[r] & \pi_1(\bar X,\bar x) \ar[r]
& \pi_1(X,\bar x)\ar[r]^{\epsilon} & {\rm Gal}(\bar{k}/k) \ar[r] &1
}}\end{thm}
\begin{proof}
Let  $\theta \in \Pi(X,\bar x)_s(\bar k)$, with image $\gamma\in {\rm Gal}(\bar{k}/k)$. 
Then by \eqref{3.3}, $\theta$ is a tensor isomorphism between  $\rho_{\bar x}$ and 
$\gamma^*\circ \rho_{\bar x}$. On the other hand, if $\pi: Y\to X$ is a Galois covering 
under finite quotient $H$ of $\pi_1(X,\bar x)$, then $\pi_*(\sO_Y,d) \in {\sf FConn}(X)$. 
Thus $\theta$ yields an isomorphism between $\sO_{\pi^{-1}(\bar x)}$ and 
$\gamma^*(\sO_{\pi^{-1}(\bar x)})$. This implies in particular that $\theta$ 
yields an automorphism of the set $\pi^{-1}(\bar x)$, that is an element in $H$. 
Since, for  another finite quotient $K$ of $\pi_1(X,\bar x)$, with $K\surj H$, the 
construction is compatible in the pro-system $\pi_1(X,\bar x)$, and  we obtain from 
$\theta$ an element in $\pi_1(X,\bar x)$. This defines the homomorphism $\kappa$. 
As $\bar{\kappa}$ is an isomorphism, $\kappa$ is an isomorphism as well. 
\end{proof}
\noindent Applying the construction in Theorem \ref{thm2.7} to the fiber functor $\rho_{\bar x}$
one obtains an \'etale cover $\pi_{\rho_{\bar x}}:X_{\rho_{\bar x}}\longrightarrow X$.
\begin{cor}\label{cor4.5}The \'etale cover  $\pi_{\rho_{\bar x}}:X_{\rho_{\bar x}}\longrightarrow X$
is the universal pro-finite \'etale covering of $X$ based at $\bar x$.
\end{cor}
\begin{proof}
Let $\tilde{\pi}_{\bar x}: \tilde X_{\bar x} \to X$ be the universal pro-finite \'etale covering based at $\bar x$. Then, $\bar x$ lifts to the $1$-element in $ \tilde{\pi}_{\bar x}^{-1}(\bar x)=\pi_1(X,\bar x) 
$ which we denote by $\tilde x$. Thus $\tilde x \in  \tilde X_{\bar x}(\bar k)$. 
By Theorem \ref{thm4.4},  $\tilde X_{\bar x}$ and $X_{\rho_{\bar x}}$ are both Galois covers  of  $X$ under the same pro-finite group via $\kappa$. Thus an isomorphism between the two covers over $X$ is uniquely determined by the image of $\tilde x$, which we determine to be the $1$-element in 
$\Pi(X, \bar x)(\bar{k})= \rho_{\bar x}^{-1}(\bar x)\subset X_{\rho_{\bar x}}(\bar{k})$. 
\end{proof}
\noindent We deduce now the main corollary of the sections 3 and 4.
\begin{cor} \label{cor4.6}
Let  $X/k$, $\bar x, \rho_{\bar x}$ be as in Theorem \ref{thm4.4}. 
Then there is a one to one correspondence between  
splittings of $\epsilon$ as pro-finite groups up to conjugation by $\pi_1(\bar{X}, \bar{x})$ and 
neutral  fiber functors ${\sf FConn}(X) \to {\sf Vec}_k$. 
\end{cor}
\begin{proof}
By Theorem \ref{thm4.4}, a splitting of $\epsilon$ of pro-finite groups 
is equivalent to a splitting  of 
$\Pi(X,\rho_x)_s (\bar k)\to  {\rm Gal}(\bar{k}/k)$. 
By Theorem \ref{thm3.2}, 2), the latter is equivalent to a splitting of 
$(t,s): \Pi(X,\rho_{\bar x}) \to \bar x\times_k \bar x$ as $k$-groupoid schemes 
acting on $\bar x$. By Theorem \ref{thm3.2}, 3), such splittings are, up to conjugation by
$\Pi(\bar X,\bar x)(\bar k)$,  in one  to one correspondence with
 neutral fiber functors on ${\sf FConn}(X)$. 

\end{proof}

\section{Applications} 

In his letter to G. Faltings  dated June 27-th, 1983 (\cite{GroFa}),  Grothendieck  
conjectures  that if $X$ is a smooth projective curve of genus $\ge 2$ defined over 
$k$ of finite type over $\Q$, then conjugacy classes of sections of $\epsilon$ in 
\eqref{4.6} are in one to one correspondence with 
rational points of $X$. Thus Corollary \ref{cor4.6} implies
\begin{thm} \label{thm5.1}
If Grothendieck's section conjecture is true, then the  set of $k$-rational points of a  smooth projective curve $X$ 
of genus $g\ge 2$ over $k$ of finite type over $\Q$  
is in bijection with the set of neutral fiber functors ${\sf FConn}(X)\to {\sf Vec}_k$. 
\end{thm}
\noindent
Theorem \ref{thm5.1} says in particular that if Grothendieck conjecture was too  
optimistic, there would be an {\it exotic} neutral  fiber functor $\rho$ on ${\sf FConn}(X)$, 
that is a fiber functor which is not {\it geometric}, i.e. not of the shape 
$\rho=\rho_x, \ \rho_x((V,\nabla))=V|_x$ for some rational point $x\in X(k)$. \\[.2cm]
On the other hand, if $\rho$ is a neutral fiber functor, then by 
Lemma \ref{lem4.3}, the Tannaka group scheme 
 $\Pi(X,\rho)$ has the property  $\Pi(X,\rho)(\bar{k})= \pi_1(X\times_k \bar{k}, x).$
The construction of Theorem \ref{thm2.7} applied to $\rho$ yields a $\Pi(X,\rho)$
principal bundle
\ga{5.1}{\pi_\rho:X_\rho\to X }
with the property
\ga{5.2}{H^0_{DR}(X_\rho) =k.}
In particular,  Theorem \ref{thm2.7}, 4) implies that 
the fiber functor $\rho$ is equivalent to the  neutral fiber functor 
\ga{}{\rho':M\mapsto H^0_{DR}(X_\rho, \pi_\rho^*(M)). \notag}
Thus one concludes
\begin{lem} \label{lem5.2}
Let $X$ be as in Theorem \ref{thm5.1}. Then if there is a finite 
connection $M$ (or equivalently a bundle in $\sC^N(X)$, see Theorem
 \ref{thm2.15}), such that for all connected \'etale coverings
$\pi: Y\to X$ which have the property that $\pi^*(M)$ is trivial, the 
field  $H^0_{DR}(Y)$ is a finite nontrivial extension of 
$k=H^0_{DR}(X)$, then ${\sf FConn}(X)$ has no neutral
 fiber functor.  In particular in this case,  $X(k)=\emptyset$. 
\end{lem}
\begin{proof} If $X$ had a $k$-rational point, say $x$, then one would have the equality
  $H_{DR}^0(X_{\rho_x})=k$.
For the finite category $\sT:=\langle M\rangle$ generated by $M$,  the canonical
morphism $X_{\rho_x}\to X_{\sT,\rho_x}$ is surjective, hence
$H^0_{DR}(X_{\sT,\rho_x})=k$, which contradicts the assumption. 
Thus $X$ has no $k$-rational point.
\end{proof}
\noindent If $\rho=\rho_x$ for some rational point $x\in X(k)$, then as discussed in Remarks \ref{rmk2.13},
\eqref{2.29} one has  $\pi_{\rho_x}^{-1}(x)=\Pi(X,\rho_x)$. In particular, $X_{\rho_x}$  carries the rational pro-point  
 $x_{\rho_x}\in X_{\rho_x}(k)$ as the unit element of $\Pi(X,\rho_x)$ with 
\ga{5.3}{\pi_{\rho_x}(x_{\rho_x})=x.}
We use here the terminology ``$k$-rational pro-point'' to denote a $k$-rational point of
 the scheme $X_{\rho_x}$, and to indicate at the same time that 
it defines a pro-system of points $X_H\in X_H(k)$ for all finite quotients $H$ of $\Pi(X, \rho_x)$.\\[.2cm]
We choose now a geometric point $\bar{x}\to x$ with residue field $\bar{k}$. 
Then by Corollary \ref{cor4.6}, \eqref{4.6} splits, and in view of Theorem \ref{thm2.11} and of
 Corollary \ref{cor4.5},
we obtain that the  universal pro-finite \'etale covering $\tilde{X}_{\bar{x}}$ based at $\bar{x}$   fulfills
\ga{5.4}{\tilde{X}_{\bar{x}}\cong X_{\rho_x}\times_k \bar{k} \ {\rm as} 
\ \bar{k}-{\rm schemes}.}
Recall that if $\bar{y}\to X$ is another geometric point with residue field 
$\bar{k}$, then there is a path $\gamma_{\bar{x}\bar{y}}
\in  \pi_1(X, \bar{x}, \bar{y})$, where $\pi_1(X, \bar x, \bar y)$ is Grothendieck's groupoid,  which yields an isomorphism 
\ga{5.5}{\gamma_{\bar{x}\bar{y}}: \tilde{X}_{\bar{x}}\cong 
\tilde{X}_{\bar{y}}\ {\rm as} \ \bar{k}-{\rm schemes}.}
 We summarize:
\begin{prop} \label{prop5.3}
Let $X/k$ be as in Notations \ref{nota4.1}. Let  $\bar{x}\to X$ be a geometric point
 with residue field $\bar{k}$. Let 
$ \tilde{X}_{\bar{x}}$ be the universal pro-finite \'etale covering of $X$ based at $\bar x$,
considered as a $\bar{k}$-scheme. Then if $\bar{x}$ descends to a $k$-rational point 
$x\in X(k)$, then $ \tilde{X}_{\bar{x}}$ descends to a $k$-scheme $X_x$, 
i.e. $\tilde{X}_{\bar{x}}\cong X_x\times_k \bar{k}$,  with 
the property that $x$ lifts to a point in $ X_x(k)$. Moreover, if 
$\bar{y}$ is another $\bar{k}$-point of $X$, choosing a path in
 $\pi_1(X,\bar{x},\bar{y})$ yields a $\bar{k}$-isomorphism 
$\tilde{X}_{\bar{x}}\cong \tilde{X}_{\bar{y}}$, and a splitting of $\epsilon$ 
for Grothendieck's sequence \eqref{4.6} based in $\bar{y}$. The splitting 
yields the $k$-structure $X_x$ on 
$ \tilde{X}_{\bar{y}}$ as well, and a $k$-rational pro-point in  $X_x(k)$ mapping to $x$.
\end{prop}
\begin{rmk} \label{rmk5.4}
We remark that Proposition \ref{prop5.3} yields  a positive answer to the conjecture 
formulated by  Grothendieck  \cite{GroFa}, page 8, once we know the existence of
 $x\in X(k)$. Loosely speaking, the $k$-linear Tannaka formalism which we introduce here in the study 
 of Grothendieck's arithmetic fundamental group, allows us to 
 easily lift $k$-rational  points to $k$-rational pro-points. ``Easily'' means here of course via 
 a careful analysis of the Tannaka formalism as developed in sections 2 and 4.
   However, it does not allow 
 to find them. For this, one needs the arithmetic of the ground field $k$ and the geometry 
 of the $X$ considered. 
\end{rmk}
\noindent
The aim of the rest of the article is to show that Proposition \ref{prop5.3} 
remains true if $X$ is an affine curve and one has a neutral fiber functor 
which is coming from a rational point at $\infty$. 
\begin{notas} \label{notas5.5}
We assume in the rest of the article that $X$ is a smooth affine curve defined 
over a field $k$ of characteristic 0, and that there is a $k$-rational point $x$ on 
the smooth compactification of $X$. We set $X'=X\cup \{x\}$. Let $\widehat{\sO}_x$ 
be the local formal ring at $x$, and $\widehat{K}_x$ be its field a fractions. 
We choose a local paramater so $ 
\widehat{\sO}_x \cong k[[t]]$ and $\widehat{K}_x\cong k((t))$.\\[.2cm]
Let $D^\times:=\Spec(k((t)))\subset D:=\Spec(k[[t]])$ be the punctured formal disk
embedded in the formal disk. 
In \cite[Section~2.1]{KaGal}, Katz defines the $k$-linear rigid tensor category  ${\sf Conn}(D^\times)$ 
of $t$-adically continuous connections on $D^\times$ as follows. Objects are pairs $M=(V,\nabla)$, where $V$ is a 
finite dimensional vector space on $k((t))$ and $\nabla: V\to V\otimes_{k((t))} \omega_{k((t))}$ a $t$-adically continuous connection, with $\omega_{k((t))}\cong k((t))dt$. 
The  Hom-Sets are flat morphisms.  We define ${\sf FConn}(D^\times)$ in the obvious way, as being the full sub-category of ${\sf Conn}(D^\times)$ 
spanned by Weil-finite objects, with definition as in  Definition \ref{defn2.3}.
One has the restriction functor
\ga{5.6}{{\sf FConn}(X) \xrightarrow{{\sf rest}} 
{\sf FConn}(D^\times).}
Let us recall the constructions of 
Deligne \cite[Section~15, 28-36]{DeP}, and Katz
\cite[Section~2.4]{KaGal}, which are depending on the $t$ chosen. There is a
 functor of $k$-linear abelian rigid tensor categories
\ga{5.7}{ {\sf FConn}(D^\times)\xrightarrow{DK} {\sf FConn}(\G_m), \ \G_m:={\rm Spec}(k[t, t^{-1}])}
which is the restriction to ${\sf FConn}$ of a functor defined on the categories of regular singular 
connections on $D^\times$ and $\G_m$ (by Deligne), and even of all connections (by Katz). It is
    characterized  by the property that it is additive, and functorial 
for inverse and direct images by ${\rm Spec}(k'((u)))\to D^\times$, $u^n=t$ where $k'\supset k$ 
is a finite field extension which contains the $n$-th rooths of unity. 
The choice of a rational point $a\in \G_m(k)$, for example $a=1$, defines a fiber functor 
$\rho_{\G_m}: {\sf FConn}(\G_m)\to {\sf Vec}_k$ by assigning $V|_a$ to $(V, \nabla)$. 
We denote the resulting functor $${\sf FConn}(D^\times)\stackrel{DK}\to
{\sf FConn}(\G_m)\stackrel {\rho_{\G_m}}\to {\sf Vec}_k$$ by $\varphi$.
The composite functor ${\sf FConn}(X)\stackrel{\sf rest}\to {\sf FConn}(D^\times)
\to {\sf Vec}_k$ will  be denoted by $\eta$:
\ga{5.8}{\xymatrix{ {\sf FConn}(X) \ar[d]_{{\sf rest}}\ar[rd]^{\eta}\\
{\sf FConn}(D^\times)\ar[r]_{\varphi}&   {\sf Vec}_k .}}
\end{notas}
\noindent
We apply the construction of Theorem \ref{thm2.7} to the 
 pair $({\sf FConn}(X),\eta)$ to obtain a morphism $ X_{\eta}\to X$.
The construction of Theorem \ref{thm2.7}, while applied to $D^\times={\rm Spec}(k((t)))$ and the 
pair $({\sf FConn}(D^\times), \varphi)$, yields 
a covering $D^\times _\varphi\to D^\times$. Indeed, as the connections are $t$-adically continuous, 
the finite sub-category spanned by a finite connection $M$ yields a finite field extension $k((t))\subset L_M$, which defines $D^\times_M:={\rm Spec}(L_M)$, and $D^\times_\varphi=\varprojlim_{M} D^\times_M$.  
The Theorem \ref{thm2.7}, 2) applies here to give $H^0_{DR}(D^\times_\varphi)=k$. 
Thus, \eqref{5.8} implies 
that
we have the following diagram
\ga{5.9}{\xymatrix{ D^\times_{{\varphi}}
 \ar[d] \ar[r] & \ar@{} [dr] |\square X_\eta \times_X D^\times \ar[d] \ar[r] &  X_\eta \ar[d] \\ D^\times \ar[r]_{=} & D^\times \ar[r] & X
}}
We now introduce the partial compactification in $x$. Since the morphism 
$D^\times_{\varphi}\to D^\times$ is indeed a pro-system of ${\rm Spec}(L_M)$ of $D^\times$,
we define its compactification as the corresponding pro-system of the normalizations of $D$ in 
$L_M$. Similarly  the compactification of
 $X_{\eta}\to X$ is the pro-system of the normalizations of $X'$ in the function fields of the coverings. We denote by $'$ the compactifications. This yields 
\ga{5.10}{ (D^\times_{\varphi})'\to D, \ (X_\eta)' \to X',}
and the diagram 
\ga{5.11}{\xymatrix{ (D^\times_{{\varphi}})'
 \ar[d] \ar[r] & \ar@{} [dr] |\square (X_\eta)' \times_{X'} D^\times \ar[d] \ar[r] &  (X_\eta)' \ar[d] \\ D \ar[r]_{=} & D \ar[r] & X'.
}}
Note the rational point  $x\in X'(k)$ lies in $D(k)$, with defining maximal ideal $\langle t \rangle$.  
 \begin{lem} \label{lem5.6}
The point $x\in D(k)$ lifts to a pro-point $x'_\varphi  \in (D^\times_\varphi)'(k)$, and thus defines a rational 
pro-point $x'_\eta$ in $(X_{\eta})'(k)$. 
\end{lem}
\begin{proof}
Etale coverings of $D$ are disjoint unions of coverings of the shape ${\rm Spec}(k'((u)))$ 
for some finite  field extension $k'\supset k$ and $u^n=t$ for some $n\in \N\setminus \{0\}$. 
Notice that $H^0_{DR}({\rm Spec}(k'((u))))=k'$.
On the other hand, Theorem \ref{thm2.7}, 2) implies that
$H^0_{DR}(D^\times_{\epsilon(D^\times)})=k$. This means that
in the pro-system defining $D^\times_{\varphi}$, there are only connected
 coverings with $k'=k$. Hence  $x \in D(k)$ lifts to $x'_\varphi \in (D^\times_{\varphi})'(k)$.
 Its image $x'_\eta$  in $(X_{\eta})'(k)$ is the required pro-point. 
\end{proof}
\noindent On the other hand, given a geometric point $\bar{y}\to X$ 
with residue field
$\bar{k}$, and given  the choice of a geometric point $\bar{x}\to X$ with residue field
 $\bar{k}$ above $x$, one has by \cite[Proposition~8.11]{DeGroth} that the two 
 ${\sf Vec}_{\bar{k}}$-valued fiber functors $\rho\times_k \bar{k}$ 
and $\rho_{\bar{y}}$ on $X\times_k \bar{k}$ are equivalent as the functor ${\rm Iso}^{\otimes}(\rho\times_k \bar{k}, \rho_{\bar{y}} )$ from ${\sf Sch}/\bar k$ to $\sf Sets$ 
is represented by a torsor  under the affine (even pro-finite) group scheme ${\rm Aut}^{\otimes}(\rho_x)$ 
over $\bar k$. Therefore  it splits. One obtains
\ga{5.12}{ \tilde{X}_{\bar{y}}\cong X_{\rho }\times_k \bar{k} \ {\rm as} \ \bar{k}-{\rm schemes}.}
We summarize:
\begin{thm} \label{thm5.7}
Let $X\subset X'/k, \eta$ be as in  Notations \ref{notas5.5}, and let 
 $\bar{y}\to X$ be a geometric point with residue field
$\bar{k}$. Then the universal pro-finite \'etale covering  
$\tilde{X}_{\bar{y}}$ based at $\bar y $ has a $k$-structure $X_\eta$, i.e. 
$\tilde{X}_{\bar{y}} \cong_{\bar{k}} X_\eta\times_k \bar{k}$, and the 
rational point $x\in X'(k)$ lifts to a $k$-rational pro-point of the  normalization $(X_\eta)'$ of $X'$ in $k(X_\eta)$.

\end{thm}
\begin{rmks} \label{rmk5.8}\begin{itemize}
 \item[1)]Theorem \ref{thm5.7} completes Proposition \ref{prop5.3} and yields  a positive answer
  to the conjecture formulated by  Grothendieck  \cite{GroFa}, page 8, for a section of 
  $\epsilon$ associated to a ``tangential vector'' (in the language of Deligne). 
\item[2)]
We have seen in remarks \ref{rmk2.13} 2) that a fiber functor 
${\sf FConn}(X)\to {\sf QCoh}(S)$ is always equivalent to a fiber functor which is 
cohomologically defined. 
In view of the approach presented in this article, we may ask for a 
generalization of Grothendieck's conjecture to more general fiber functors:
what are the geometric assumptions on a smooth scheme $X$ of finite type over $k$  
and the arithmetic assumption on $k$ which force the following to be true:
there is an $S$-valued fiber functor of ${\sf FConn}(X)$ if and only if there is an
$S$-point of $X$? 
One could go on and ask for pro-points as well.
A generalization would be saying that if  an  $S$-valued fiber functor $\rho$ 
comes from a $S$-point of $X$, then this $S$-point lifts as an $S$-pro-point of $X_\rho$ 
defined in Theorem \ref{thm2.7}. \end{itemize}
\end{rmks}

\bibliographystyle{plain}
\renewcommand\refname{References}

\end{document}